\documentclass[11pt]{article}
\usepackage[utf8]{inputenc}
\usepackage[T1]{fontenc}
\usepackage[colorlinks=true]{hyperref}
\hypersetup{
 urlcolor=blue,
 citecolor=blue,
 linkcolor=black
}
\usepackage[dvips]{graphicx}
\usepackage{slashbox}
\usepackage{textcomp}
\usepackage{amsmath, amsbsy, amsfonts, amssymb}
\usepackage{latexsym}
\usepackage[mathscr]{eucal}
\usepackage{mathtools}
\mathtoolsset{showonlyrefs=false} 
\usepackage[usenames]{color}
\usepackage{subcaption}
\usepackage[dvipsnames,table]{xcolor}
\usepackage{url}
\usepackage{tikz}

\usepackage{mathabx}

\usepackage{algorithm,algorithmicx,algpseudocode}

\newcommand{\ac}[1]{\textcolor{black}{#1}}
\newcommand{\dw}[1]{\textcolor{black}{#1}}
\newcommand{\dwo}[1]{\textcolor{black}{#1}}
\newcommand{\jn}[1]{\textcolor{black}{#1}}
\newcommand{\om}[1]{{\color{black}{#1}}}
\def\anew{\color{black} }
\usepackage{fullpage}

\begin{document}
\bibliographystyle{plain}
\title
{
Nonlinear reduced models for state and parameter estimation}
\author{ 
Albert Cohen, Wolfgang Dahmen, Olga Mula and James Nichols\thanks{%
This research was supported by the ERC Adv grant BREAD; the Emergence Project of the Paris City Council ``Models and Measures''; the NSF Grants DMS ID 1720297, DMS ID 2012469, and by the SmartState and Williams-Hedberg Foundation.
}  }

\hbadness=10000
\vbadness=10000
\newtheorem{lemma}{Lemma}[section]
\newtheorem{prop}[lemma]{Proposition}
\newtheorem{cor}[lemma]{Corollary}
\newtheorem{theorem}[lemma]{Theorem}
\newtheorem{remark}[lemma]{Remark}
\newtheorem{example}[lemma]{Example}
\newtheorem{definition}[lemma]{Definition}
\newtheorem{proper}[lemma]{Properties}
\newtheorem{assumption}[lemma]{Assumption}
%
\newenvironment{disarray}{\everymath{\displaystyle\everymath{}}\array}{\endarray}

\def\vp{\varphi}
\def\<{\langle}
\def\>{\rangle}
\def\t{\tilde}
\def\i{\infty}
\def\e{\varepsilon}
\def\sm{\setminus}
\def\nl{\newline}
\def\o{\overline}
\def\wt{\widetilde}
\def\wh{\widehat}
\def\cT{{\cal T}}
\def\cA{{\cal A}}
\def\cI{{\cal I}}
\def\cV{{\cal V}}
\def\cB{{\cal B}}
\def\cF{{\cal F}}

\def\cR{{\cal R}}
\def\cD{{\cal D}}
\def\cP{{\cal P}}
\def\cJ{{\cal J}}
\def\cM{{\cal M}}
\def\cO{{\cal O}}
\def\Chi{\raise .3ex
\hbox{\large $\chi$}} \def\vp{\varphi}
\def\lsima{\hbox{\kern -.6em\raisebox{-1ex}{$~\stackrel{\textstyle<}{\sim}~$}}\kern -.4em}
\def\lsim{\hbox{\kern -.2em\raisebox{-1ex}{$~\stackrel{\textstyle<}{\sim}~$}}\kern -.2em}
\def\[{\Bigl [}
\def\]{\Bigr ]}
\def\({\Bigl (}
\def\){\Bigr )}
\def\[{\Bigl [}
\def\]{\Bigr ]}
\def\({\Bigl (}
\def\){\Bigr )}
\def\L{\pounds}
\def\pr{{\rm Prob}}
\newcommand{\cs}[1]{{\color{magenta}{#1}}}
\def\ds{\displaystyle}
\def\ev#1{\vec{#1}}     
\newcommand{\lt}{\ell^{2}(\nabla)}
\def\Supp#1{{\rm supp\,}{#1}}
\def\R{\mathbb{R}}
\def\E{\mathbb{E}}
\def\nl{\newline}
\def\T{{\relax\ifmmode I\!\!\hspace{-1pt}T\else$I\!\!\hspace{-1pt}T$\fi}}
\def\N{\mathbb{N}}
\def\Z{\mathbb{Z}}
\def\N{\mathbb{N}}
\def\Zd{\Z^d}
\def\Q{\mathbb{Q}}
\def\C{\mathbb{C}}
\def\Rd{\R^d}
\def\gsim{\mathrel{\raisebox{-4pt}{$\stackrel{\textstyle>}{\sim}$}}}
\def\sime{\raisebox{0ex}{$~\stackrel{\textstyle\sim}{=}~$}}
\def\lsim{\raisebox{-1ex}{$~\stackrel{\textstyle<}{\sim}~$}}
\def\div{\mbox{ div }}
\def\M{M}  \def\NN{N}                  
\def\L{{\ell}}               
\def\Le{{\ell^1}}            
\def\Lz{{\ell^2}}
\def\Let{{\tilde\ell^1}}     
\def\Lzt{{\tilde\ell^2}}
\def\Ltw{\ell^\tau^w(\nabla)}
\def\t#1{\tilde{#1}}
\def\la{\lambda}
\def\La{\Lambda}
\def\ga{\gamma}
\def\BV{{\rm BV}}
\def\Ga{\eta}
\def\al{\alpha}
\def\cZ{{\cal Z}}
\def\cA{{\cal A}}
\def\cU{{\cal U}}
\def\ms{{\rm ms}}
\def\wc{{\rm wc}}
\def\argmin{\mathop{\rm argmin}}
\def\argmax{\mathop{\rm argmax}}
\def\prob{\mathop{\rm prob}}
\def\A{\mathop{\rm Alg}}

\def \bphi{{\bf\phi}}

\def\cO{{\cal O}}
\def\cA{{\cal A}}
\def\cC{{\cal C}}
\def\cS{{\cal F}}
\def\bu{{\bf u}}
\def\bz{{\bf z}}
\def\bZ{{\bf Z}}
\def\bI{{\bf I}}
\def\cE{{\cal E}}
\def\cD{{\cal D}}
\def\cG{{\cal G}}
\def\cI{{\cal I}}
\def\cJ{{\cal J}}
\def\cM{{\cal M}}
\def\cN{{\cal N}}
\def\cT{{\cal T}}
\def\cU{{\cal U}}
\def\cV{{\cal V}}
\def\cW{{\cal W}}
\def\cL{{\cal L}}
\def\cB{{\cal B}}
\def\cG{{\cal G}}
\def\cK{{\cal K}}
\def\cS{{\cal S}}
\def\cP{{\cal P}}
\def\cQ{{\cal Q}}
\def\cR{{\cal R}}
\def\cU{{\cal U}}
\def\bL{{\bf L}}
\def\bl{{\bf l}}
\def\bK{{\bf K}}
\def\bC{{\bf C}}
\def\X{X\in\{L,R\}}
\def\ph{{\varphi}}
\def\D{{\Delta}}
\def\H{{\cal H}}
\def\bM{{\bf M}}
\def\bx{{\bf x}}
\def\bj{{\bf j}}
\def\bG{{\bf G}}
\def\bQ{{\bf Q}}
\def\bS{{\bf S}}
\def\bP{{\bf P}}
\def\bW{{\bf W}}
\def\bX{{\bf X}}
\def\bT{{\bf T}}
\def\bV{{\bf V}}
\def\bv{{\bf v}}
\def\bt{{\bf t}}
\def\bz{{\bf z}}
\def\bw{{\bf w}}
\def \meas {{\rm meas}}
\def\rhom{{\rho^m}}
\def\diff{\hbox{\tiny $\Delta$}}
\def\EE{{\rm Exp}}
\def\lll{\langle}
\def\argmin{\mathop{\rm argmin}}
\def\argmax{\mathop{\rm argmax}}
\def\dJ{\nabla}
\newcommand{\ba}{{\bf a}}
\newcommand{\bb}{{\bf b}}
\newcommand{\bc}{{\bf c}}
\newcommand{\bd}{{\bf d}}
\newcommand{\bs}{{\bf s}}
\newcommand{\bff}{{\bf f}}
\newcommand{\bp}{{\bf p}}
\newcommand{\bg}{{\bf g}}
\newcommand{\by}{{\bf y}}
\newcommand{\br}{{\bf r}}
\newcommand{\be}{\begin{equation}}
\newcommand{\ee}{\end{equation}}
\newcommand{\bea}{$$ \begin{disarray}{lll}}
\newcommand{\eea}{\end{disarray} $$}
\def \Vol{\mathop{\rm  Vol}}
\def \mes{\mathop{\rm mes}}
\def\rad{\mathop{\rm rad}}
\def \Prob{\mathop{\rm  Prob}}
\def \exp{\mathop{\rm    exp}}
\def \sign{\mathop{\rm   sign}}
\newcommand{\mult}{\mathop{\rm   mult}}
\newcommand{\one}{\mathop{\rm   one}}

\def \sp{\mathop{\rm   span}}
\def \vphi{{\varphi}}
\def \csp{\overline \mathop{\rm   span}}

\newcommand{\miny}{\textsc{MIN-Y}}
\newcommand{\minv}{\textsc{MIN-V}}
\newcommand{\altmin}{Alt-MIN}

%
%
\newcommand{\beqn}{\begin{equation}}
\newcommand{\eeqn}{\end{equation}}
\def\beginproof{\noindent{\bf Proof:}~ }
\def\endproof{\hfill\rule{1.5mm}{1.5mm}\\[2mm]}

\newcommand{\utr}{u^{\rm true}}
\newcommand{\Cor}{\kappa}

\newenvironment{Proof}{\noindent{\bf Proof:}\quad}{\endproof}

\renewcommand{\theequation}{\thesection.\arabic{equation}}
\renewcommand{\thefigure}{\thesection.\arabic{figure}}

\makeatletter
\@addtoreset{equation}{section}
\makeatother

\newcommand\abs[1]{\left|#1\right|}
\newcommand\clos{\mathop{\rm clos}\nolimits}
\newcommand\trunc{\mathop{\rm trunc}\nolimits}
\renewcommand\d{d}
\newcommand\dd{\mathrm d}
\newcommand\diag{\mathop{\rm diag}}
\newcommand\dist{\mathop{\rm dist}}
\newcommand\diam{\mathop{\rm diam}}
\newcommand\cond{\mathop{\rm cond}\nolimits}
\newcommand\eref[1]{{\rm (\ref{#1})}}
\newcommand{\iref}[1]{{\rm (\ref{#1})}}
\newcommand\Hnorm[1]{\norm{#1}_{H^s([0,1])}}
\def\int{\intop\limits}
\renewcommand\labelenumi{(\roman{enumi})}
\newcommand\lnorm[1]{\norm{#1}_{\ell^2(\Z)}}
\newcommand\Lnorm[1]{\norm{#1}_{L_2([0,1])}}
\newcommand\LR{{L_2(\R)}}
\newcommand\LRnorm[1]{\norm{#1}_\LR}
\newcommand\Matrix[2]{\hphantom{#1}_#2#1}
\newcommand\norm[1]{\left\|#1\right\|}
\newcommand\ogauss[1]{\left\lceil#1\right\rceil}
\newcommand{\QED}{\hfill
\raisebox{-2pt}{\rule{5.6pt}{8pt}\rule{4pt}{0pt}}%
  \smallskip\par}
\newcommand\Rscalar[1]{\scalar{#1}_\R}
\newcommand\scalar[1]{\left(#1\right)}
\newcommand\Scalar[1]{\scalar{#1}_{[0,1]}}
\newcommand\Span{\mathop{\rm span}}
\newcommand\supp{\mathop{\rm supp}}
\newcommand\ugauss[1]{\left\lfloor#1\right\rfloor}
\newcommand\with{\, : \,}
\newcommand\Null{{\bf 0}}
\newcommand\bA{{\bf A}}
\newcommand\bB{{\bf B}}
\newcommand\bR{{\bf R}}
\newcommand\bD{{\bf D}}
\newcommand\bE{{\bf E}}
\newcommand\bF{{\bf F}}
\newcommand\bH{{\bf H}}
\newcommand\bU{{\bf U}}
\newcommand\cH{{\cal H}}
\newcommand\sinc{{\rm sinc}}
\def\enorm#1{| \! | \! | #1 | \! | \! |}

\newcommand{\dm}{\frac{d-1}{d}}

\let\bm\bf
\newcommand{\balpha}{{\mbox{\boldmath$\alpha$}}}
\newcommand{\bbeta}{{\mbox{\boldmath$\beta$}}}
\newcommand{\bal}{{\mbox{\boldmath$\alpha$}}}
\newcommand{\bbi}{{\bm i}}

\def\nnew{\color{black}}
\def\mnew{\color{Blue}}

\newcommand{\dI}{\Delta}
%
%

\newcommand\rank{\mathop{\rm rank}}
\newcommand\tr{\mathop{\rm trace}}
\newcommand\ra{\mathop{\rm range}}
\newcommand\bPhi{{\bf \Phi}}
\newcommand\bPsi{{\bf \Psi}}
\newcommand\bSigma{{\bf \Sigma}}
\newcommand\bGamma{{\bf \Gamma}}
\newcommand\bbP{\mathbb{P}}
\newcommand\rd{\mathrm{d}}

 \newcommand{\mnote}[1]{\marginpar{\scriptsize \textcolor{magenta}{WD: #1}}}
 \newcommand{\jmnote}[1]{\marginpar{\scriptsize \textcolor{LimeGreen}{JN: #1}}}
  \newcommand{\omnote}[1]{\marginpar{\scriptsize \textcolor{blue}{OM: #1}}}

\maketitle
\date{}

\begin{abstract}
State estimation aims at approximately reconstructing the solution $u$ to a parametrized partial differential equation 
from $m$ linear measurements, when the parameter vector $y$ is unknown. Fast numerical recovery methods have been 
proposed in \cite{MPPY} based on reduced models which are linear spaces of moderate dimension 
$n$ which are tailored to approximate the solution manifold $\cM$ where the solution sits. These methods can be viewed as 
deterministic counterparts to Bayesian estimation approaches, and are proved to be optimal when the prior is
expressed by approximability of the solution with respect to the reduced model \cite{BCDDPW2}.
However, they are inherently limited by their linear nature, which bounds from
below their best possible performance by the Kolmogorov width $d_m(\cM)$ of the solution manifold.
In this paper we propose to break this barrier by using simple nonlinear reduced models that consist
of a finite union of linear spaces $V_k$, each having dimension at most $m$
and leading to different estimators $u_k^*$. A model selection mechanism 
based on minimizing the PDE residual over the parameter space is used to select from this collection the
final estimator $u^*$. Our analysis shows that $u^*$ meets optimal 
recovery benchmarks that are inherent to the solution manifold and not tied to its Kolmogorov width.
The residual minimization procedure is computationally simple in the relevant case of affine parameter
dependence in the PDE. In addition, it results in an estimator $y^*$ for the unknown parameter vector.
In this setting, we also discuss an alternating minimization \ac{(coordinate descent)} algorithm for joint state and parameter
estimation, that potentially improves the quality of both estimators. 
\end{abstract}
 
\section{Introduction}

\subsection{Parametrized PDEs and inverse problems}

Parametrized partial differential equations are of common used to model
complex physical systems. Such equations can generally be written in abstract form as
\be
\cP(u,y)=0,
\label{genpar}
\ee
where $y=(y_1,\dots,y_d)$ is a vector of
scalar parameters ranging in some domain $Y\subset \R^d$.
We assume well-posedness, that is, for any $y\in Y$  the problem admits a unique solution $u=u(y)$ in some 
Hilbert space $V$. 
We may therefore consider the {\it parameter to solution map}
\be
y\mapsto u(y),
\label{solmap}
\ee
from $Y$ to $V$, which is typically nonlinear, as well as the {\it solution manifold}
\be
\cM:=\{u(y) \, : \, y\in Y\}\subset V
\label{solman}
\ee
that describes the collection of all admissible solutions. Throughout this paper,
we assume that $Y$ is compact in $\R^d$ and that the map \iref{solmap}
is continuous. Therefore $\cM$ is a compact set of $V$. 
We sometimes refer to the solution $u(y)$ as the {\it state} 
of the system for the given parameter vector $y$.

The parameters are used to represent physical quantities such as 
diffusivity, viscosity, velocity, source terms, or the geometry of the physical domain
in which the PDE is posed.  In several relevant instances, $y$ may be high or even
countably infinite dimensional, that is, $d\gg1$ or $d=\infty$. 

In this paper, we are interested in {\it inverse problems} which occur when only a  vector of {\it linear} measurements
\be
z=(z_1,\dots,z_m)\in \R^m, \quad z_i=\ell_i(u),\quad i=1,\dots,m,
\ee
is observed, where each $\ell_i\in V'$ is a known continuous linear functional on $V$. \ac{We also sometimes use the notation
\be
z=\ell(u), \quad \ell=(\ell_1,\dots,\ell_m).
\ee
}
One wishes to recover from $z$ the unknown state $u\in \cM$ or even the underlying parameter vector $y\in Y$ for which $u=u(y)$.
Therefore, in an idealized setting, one partially observes the result
of the composition map
\be
y\in Y \mapsto u \in \cM \mapsto z\in \R^m.
\ee
for the unknown $y$. More realistically, the measurements may be affected by additive noise 
\be
z_i=\ell_i(u)+\eta_i,
\ee
and the model itself might be biased, meaning that the true state $u$ deviates from the solution manifold $\cM$
by some amount. Thus, two types of inverse problems may be considered:
\begin{enumerate}
\item
State estimation: recover an approximation $u^*$ of \ac{the state} $u$ from the observation \ac{$z=\ell(u)$}. This is a linear inverse problem,
in which the prior information on $u$ is given by the manifold $\cM$ which has a complex geometry and 
is not explicitly known.
\item
Parameter estimation: recover an approximation $y^*$ of the parameter $y$ from the observation \ac{$z = \ell(u)$ when $u=u(y)$}. This is a nonlinear inverse problem, for which the prior information available on $y$ is given by the domain $Y$.
\end{enumerate}

These problems become severely ill-posed when $\dw{Y}$ has dimension $d>m$. For this reason,
they are often addressed through Bayesian approaches \cite{DS,St}: a prior probability distribution $P_y$ being assumed on $y\in Y$
(thus inducing a push forward distribution $P_u$ for $u\in \cM$), the objective is to understand the {\it posterior} distributions
of $y$ or $u$ conditioned by the observations $z$ in order to compute plausible solutions $y^*$ or $u^*$ under such probabilistic priors.
The accuracy of these solutions should therefore be assessed in some average sense.

In this paper, we do not follow this avenue: the only priors made on $y$ and $u$ are their
membership to $Y$ and $\cM$. We are interested in developping practical estimation methods that offer
uniform recovery guarantees under such deterministic priors in the form of upper bounds on the worst case
error for the estimators over all $y\in Y$ or $u\in \cM$. \ac{ We also aim to understand whether our error
bounds are optimal in some sense. Our primary focus will actually be on state estimation (i).
Nevertheless we present in \S\ref{sec:parest} several implications on parameter estimation (ii), which to our
knowledge are new.  For state estimation, error bounds have recently been established  
for a class of methods based on {\it linear reduced modeling}, as we  recall next.}

\subsection{Reduced models: the PBDW method}

In several relevant instances, the particular parametrized PDE structure allows one  to construct linear
spaces $V_n$ of moderate dimension $n$ that are specifically tailored to the approximation 
of the solution manifold $\cM$, in the sense that
\be
{\rm dist}(\cM,V_n)=\max_{u\in\cM} \min_{v\in V_n} \|u-v\| \leq  \e_n,
\ee
where $\e_n$ is a certified bound that decays with $n$ significanly faster than when using for $V_n$ classical approximation
spaces  {of dimension $n$} such as finite elements, algebraic or trigonometric polynomials, \dw{or} spline functions. Throughout this paper
\be
\dw{\langle \cdot,\cdot\rangle_V^{1/2} =:} \|\cdot\|=\|\cdot\|_V,
\ee
denotes the norm of the Hilbert space $V$. The natural benchmark for such approximation spaces is the
Kolmogorov $n$-width
\be
d_n(\cM):=\min_{\dim(E)=n} {\rm dist}(\cM,E).
\ee
The space $E_n$ that achieves the above minimum is thus the best possible reduced model  {for approximating all of $\cM$},
however it is computationally out of reach.

One instance of computational reduced model spaces is generated by {\it sparse polynomial approximations}
of the form
\be
u_n(y)=\sum_{\nu\in \Lambda_n} u_\nu y^\nu, \quad y^\nu:=\prod_{j\geq 1}y_j^{\nu_j},
\ee
where $\Lambda_n$ is a conveniently chosen set of multi-indices such that $\#(\Lambda_n)=n$.
Such approximations can be derived,  for example,  by best $n$-term truncations of infinite Taylor or 
orthogonal polynomial expansions. We refer to \cite{CD,CDS1} where convergence estimates of the form
\be
\sup_{y\in Y}\|u(y)-u_n(y)\|\leq Cn^{-s},
\ee
are established for some $s>0$ even when $d=\infty$. Therefore, the space $V_n:={\rm span}\{u_\nu\, : \, \nu\in\Lambda_n\}$
approximates the solution manifold with accuracy $\e_n=Cn^{-s}$. 

Another instance, known as {\it reduced basis approximation}, consists of using spaces
of the form 
\be
V_n:={\rm span}\{u_1,\dots,u_n\},
\ee 
where $u_i=u(y^i)\in\cM$ are instances
of solutions corresponding to a particular selection of parameter values $y^i\in Y$ (see \cite{MPT,RHP,S}). One typical selection procedure is based on a greedy algorithm:
 {one picks}  $y^k$ such that  {$u_k= u(y^k)$} is furthest away from the previously constructed space $V_{k-1}$ in the sense of maximizing a computable 
and tight a-posteriori bound of the projection error $\|u(y)-P_{V_{k-1}}u(y)\|$ over
a sufficiently fine discrete training set $\t Y\subset Y$. In turn, this method also delivers
a computable upper estimate $\e_k$ for ${\rm dist}(\cM,V_k)$. It was proved in 
\cite{BCDDPW1,DPW2} that the reduced basis spaces resulting from this greedy algorithm have near-optimal approximation property, in the sense that if $d_n(\cM)$ has a certain polynomial or exponential
rate of decay as $n\to \infty$, then the same rate is achieved by ${\rm dist}(\cM,V_n)$.

In both cases, these reduced models come in the form of a hierarchy $(V_n)_{n\geq 1}$, with computable decreasing
error bounds $(\e_n)_{n\geq 1}$, where $n$ corresponds to the level of truncation 
in the first case and the step of the greedy algorithm in the second case. Given a reduced model $V_n$, 
one way of tackling the state estimation problem is to replace the complex solution manifold $\cM$
by the simpler prior class described by the cylinder
\be
\cK=\cK(V_n,\e_n)=\{v\in V\; : \; \dist(v,V_n)\leq \e_n\}.
\ee
that contains $\cM$. The set $\cK$ therefore reflects the approximability of $\cM$ by $V_n$.  
This point of view leads to the {\it Parametrized Background Data Weak} (PBDW)
method introduced in \cite{MPPY}, also termed as {\it one space method} and further analyzed in \cite{BCDDPW2}, that we recall 
below in a nutshell.

In the noiseless case, the knowledge of $z=(z_i)_{i=1,\dots,m}$ is equivalent 
to that of the orthogonal projection $w=P_{W} u$, where
\be
W:={\rm span}\{\omega_1,\dots,\omega_m\}
\ee
and $\omega_i\in V$ are the Riesz representers of the linear functionals $\ell_i$, that is
\be
\ell_i(v)=\<\omega_i,v\>, \quad v\in V.
\ee
Thus, the data indicates that $u$ belongs to the affine space 
\be
V_w:=w+W^\perp.
\ee
Combining this information with the prior class $\cK$, the unknown state thus belongs 
to the ellipsoid
\be
\cK_w:=\cK\cap V_w=\{v\in \cK\; : \; P_Wv=w\}.
\ee
For this posterior class $\cK_w$, the {\it optimal recovery} estimator $u^*$ that 
minimizes the worst case error $\max_{u\in \cK_w}\|u-u^*\|$
is therefore the center of the ellipsoid, which is equivalently given by 
\be 
\label{u_star_def}
u^*=u^*(w):={\rm argmin}\{ \|v-P_{V_n} v\| \; : \; P_{W}  v=w\}.
\ee
It can be computed from the data $w$ in an elementary manner by solving a finite set of linear
equations. The worst case performance for this estimator, both over $\cK$ and $\cK_w$,
for any $w$, is thus given by the half-diameter of the ellipsoid which is the product of the
width $\e_n$ of $\cK$ \dw{and} the quantity
\be
\mu_n=\mu(V_n,W):=\max_{v\in V_n} \frac {\|v\|}{\|P_W v\|}.
\ee
\dw{Note that $\mu_n$} \om{is the inverse of the cosine of the angle between $V_n$ and $W$.}  For $n\ge 1$,
this quantity can be computed as the inverse of the smallest singular value
of the $n\times m$ cross-Gramian matrix with entries $\<\phi_i,\psi_j\>$ between any pair of
orthonormal bases $(\phi_i)_{i=1,\dots,n}$ and $(\psi_j)_{j=1,\dots,m}$ of $V_n$ and $W$, respectively.
\ac{It is readily seen that one also has
\be
\mu_n=\max_{w\in W^\perp} \frac {\|w\|}{\|P_{V_n^\perp} w\|},
\ee
allowing us to extend the above definition to the case of the zero-dimensional space $V_n = \{0\}$
for which $\mu(\{0\},W) =1$.} 

Since $\cM\subset \cK$, the worst case error bound over $\cM$ of the estimator, defined as
\be
E_{wc}:=\max_{u\in \cM} \|u-u^*(P_Wu)\|,
\ee
satisfies the error bound
\be 
\label{rec_error}
E_{wc}=\max_{u\in \cM} \|u-u^*(P_Wu)\|  \leq\max_{u\in \cK} \|u-u^*(P_Wu)\| \ac{= \mu_n \e_n}.
\ee

\begin{remark} 
\label{rembiasnoise}
The estimation map $w\mapsto u^*(w)$ is linear with norm $\mu_n$
and does not depend on $\e_n$. It thus satisfies, for any individual $u\in V$ and $\eta\in W$,
\be
\|u-u^*(P_Wu+\eta)\| \leq \mu_n({\rm dist}(u,V_n)+\|\eta\|),
\ee
We may therefore account for an additional measurement noise 
and model bias: if the observation is $w=P_Wu+\eta$ with $\|\eta\|\leq \e_{noise}$,
and if the true states do not lie in $\cM$ but satisfy ${\rm dist}(u,\cM)\leq \e_{model}$,
the guaranteed error bound \iref{rec_error} should be modified into
\be
\|u-u^*(w)\| \leq \mu_n(\e_n+\e_{noise}+\e_{model}).
\ee
\ac{In practice, the noise component $\eta\in W$ typically results from 
a noise vector $\o \eta\in\R^m$ affecting the
observation $z$ according to $z= \ell(u) +\o \eta$.
Assuming a bound $\|\o \eta\|_2 \leq \o \e_{noise}$ where
$\|\cdot\|_2$ is the Euclidean norm in $\R^m$, we thus receive the
above error bound with $\e_{noise}:=\|M\| \o\e_{noise}$,
where $M\in \R^{m\times m}$ is the matrix that transforms the representer basis $\omega = \{\omega_1,\ldots,\omega_m\}$
into an orthonormal basis $\psi = \{\psi_1,\ldots,\psi_m\}$ of $W$.}
\dw{Here estimation accuracy benefits
from decreasing noise without increasing computational cost. This
is in contrast to Bayesian methods for which small noise level induces
computational difficulties due to the concentration of the
posterior distribution.} \end{remark}

\begin{remark}
\label{spacediscretization}
 To bring out the essential mechanisms, we have idealized (and technically simplified)  the description of the PBDW method by omitting 
certain discretization aspects that are unavoidable in computational practice
and should be accounted for. To start with, the snapshots $u_i$ (or the polynomial
coefficients $u_\nu$) that span the reduced basis spaces $V_n$ cannot be computed
exactly, but only up to some tolerance by a numerical solver. One typical instance
is the finite element method, which yields an approximate parameter to solution map
\be
y\mapsto u_h(y) \in V_h,
\ee
where $V_h$ is a reference finite element space ensuring a prescribed accuracy
\be
\|u(y)-u_h(y)\|\leq \e_h, \quad y\in Y.
\ee
The computable states are therefore elements of the perturbed manifold
\be
\cM_h:=\{u_h(y) \, : \, y\in Y\}.
\ee
The reduced model spaces $V_n$ are low dimensional subspaces of $V_h$, and 
with certified accuracy 
\be
{\rm dist}(\cM_h,V_n)\leq \e_n.
\ee
The true states do not belong to $\cM_h$ and this deviation can therefore be interpreted
as a model bias in the sense of the previous remark with $\e_{model}=\e_h$. The
application of the PDBW also requires the introduction of the Riesz lifts $\omega_i$ in order
to define the measurement space $W$. Since we operate in the space $V_h$, these can be defined
as elements of this space satisfying
\be
\<\omega_i,v\>_{\dw{V}}=\ell_i(v), \quad v\in V_h,
\ee
thus resulting in a measurement space $W\subset V_h$. For example, if $V$ is the Sobolev spaces $H^1_0(\Omega)$
for some domain $\Omega$ and $V_h$ is a finite element subspace, the Riesz lifts are the unique solutions to the
Galerkin problem
\be
\int_{\Omega} \nabla \omega_i\nabla v=\ell_i(v), \quad v\in V_h,
\ee
and can be identified by solving $n_h\times n_h$ linear systems.  {Measuring accuracy in $V$, i.e., in a metric dictated by the 
continuous PDE model, the   idealization, to be largely maintained in what follows,  also helps understanding how to properly
adapt the background-discretization $V_h$ to the overall achievable estimation accuracy.} Other computational issues involving the 
space $V_h$ will be discussed in \S 3.4. 
\end{remark}

Note that $\mu_n\geq 1$ increases with $n$ and that its finiteness 
imposes that $\dim(V_n)\leq \dim(W)$, that is $m\geq n$. Therefore, one natural way to
decide which space $V_n$ to use is to take the value of $n\in \{0,\dots,\om{m}\}$ that minimizes
the bound $\mu_n\e_n$. This choice is somehow crude since it might not be the value
of $n$ that minimizes the true reconstruction error for a given $u\in \cM$, and for this reason 
it was referred to as a {\it poor man algorithm} in \cite{BCDDPW1}.

The PBDW approach to state estimation can be improved in various 
ways:
\begin{itemize}
\item
One variant that is relevant to the present work is to use reduced models of affine form
\be
V_n=\o u_n+\o V_n,
\ee
where $\o V_n$ is a linear space and $\o u$ is a given offset. The optimal recovery estimator 
is again defined by the minimization property \iref{u_star_def}. Its computation
amounts to the same type of linear systems and the reconstruction map $w\mapsto u^*(w)$ is now affine. The error bound
\iref{rec_error} remains valid with $\mu_n=\mu(\o V_n,W)$
and $\e_n$ a bound for ${\rm dist}(\cM,V_n)$.   
Note that $\e_n$ is also a bound for the distance of $\cM$ to the linear space $\o V_{n+1}:=\o V_n\,\oplus\, \R \o u_n$ of
dimension $n+1$. However, using instead this linear space, could result in a stability constant $\mu_{n+1}=\mu(\o V_{n+1},W)$
that is much larger than $\mu_n$, in particular, when the offset $\o u_n$ is close to $W^\perp$. 
\item
Another variant proposed in \cite{CDDFMN} consists in using a large set 
$\cT_N=\{u_i=u(y^i) \,: \, i=1,\dots,N\}$ of precomputed solutions in order to train the reconstruction maps $w\mapsto u^*(w)$
by minimizing the least-square fit $\sum_{i=1}^N \|u_i-u^*(P_Wu_i)\|^2$ over all linear or affine map,
which amounts to optimizing the choice of the space $V_n$ in the PBDW method. 
\item
Conversely, for a given reduced basis space $V_n$, it is also possible to optimize the choice of
linear functionals $(\ell_1,\dots,\ell_m)$ giving rise to the data, among a dictionary $\cD$ that represent a set of 
admissible measurement devices. The objective is to minimize the stability constant $\mu(V_n,W)$ for the resulting
space $W$, see in particular \cite{BCMN} where a greedy algorithm is proposed for selecting the $\ell_i$.
We do not take this view in the present paper and think of the space $W$ as fixed once and for all:
the measurement devices are given to us and cannot be modified.
\end{itemize}

\subsection{Objective and outline} \label{ssec:obj}

The simplicity of the PBDW method and its above variants come together with a 
fundamental limitation of its performance: since the map $w\mapsto u^*(w)$ is linear or affine,
the reconstruction necessarily belongs to an $m$ or $m+1$ dimensional space, and thefore the worst case performance
is necessarily bounded from below by the Kolmogorov width $d_m(\cM)$ or $d_{m+1}(\cM)$.
In view of this limitation, one principle objective of the present work
is to develop {\it nonlinear} state estimation techniques which provably overcome
the bottleneck of the Kolmogorov width $d_m(\cM)$. 

In \S 2, we introduce various benchmark quantities that describe the best 
possible performance of a recovery map in a worst case sense.
We first consider an idealized setting where the state $u$ is assumed to 
exactly satisfy the theoretical model described by the parametric PDE,
that is $u\in\cM$. Then we introduce similar benchmarks quantities in the presence
of model bias and measurement noise. All these quantities can be substantially
smaller than $d_m(\cM)$.

In \S 3, we discuss a nonlinear recovery method, based on
a family of affine reduced models $(V_k)_{k=1,\dots,K}$,
where each $V_k$ has dimension $n_k\leq m$ and serves 
as a local approximations to a portion $\cM_k$ of the solution manifold. 
Applying the PBDW method with each such space, results in a collection 
of state estimators $u^*_k$. The value $k$ for which the true state $u$ belongs
to $\cM_k$ being unknown, we introduce a {\em model selection} procedure in order 
to pick a value $k^*$, and define the resulting estimator $u^*=u^*_{k^*}$.
We show that this estimator has performance comparable to the 
benchmark introduced in \S 2. Such performances cannot be achieved 
by the standard PBDW method due to the above described limitations.

Model selection is a classical topic of mathematical statistics \cite{Mas},
with representative techniques such as complexity penalization or cross-validation
in which the data are used to select a proper model. Our approach differs
from these techniques in that it exploits  {(in the spirit of {\em data assimilation})} the PDE model which is available
to us, by evaluating the distance to the manifold
\be
\dist(v,\cM)=\min_{y\in Y} \|v-u(y)\|,
\label{distM}
\ee
of the different estimators $v=u^*_k$ for $k=1,\dots,K$, and picking the value
$k^*$ that minimizes it. In practice, the quantity \iref{distM}
cannot be exactly computed and we instead rely on a computable surrogate quantity 
${\cal S}(v,\cM)$ expressed in terms of the residual to the PDE, {see \S~\ref{ssec:residual}}. 

One typical instance 
where such a surrogate is available is when \iref{genpar} 
has the form of a linear operator equation
\be
A(y)u=f(y),
\ee
where $A(y)$ is boundedly invertible from $V$ to $V'$, or more generally, from $V\to Z'$ for a test space $Z$ different from $V$, uniformly over $y\in Y$. Then ${\cal S}(v,\cM)$ is obtained
by minimizing the residual 
\be
\cR(v,y)=\|A(y)v\dw{-}f(y)\|_{Z'},
\ee
over $y\in Y$. This task itself is greatly facilitated in the case where the operators $A(y)$ and
source terms $f(y)$ have affine dependence in $y$. One relevant example that has been 
often considered in the literature is the second order elliptic diffusion equation with affine diffusion coefficient,
\be
-{\rm div}(a\nabla u)=f(y), \quad a=a(y)=\o a+\sum_{j=1}^d y_j\psi_j.
\label{ellip}
\ee

In \S 4, we discuss the more direct approach for both state and parameter estimation
based on minimizing $\cR(v,y)$ over both $y\in Y$ and $v\in w+W^\perp$.
The associated alternating minimization algorithm amounts 
to a simple succession of quadratic problems 
in the particular case of linear PDE's with affine parameter dependence. Such an algorithm is not guaranteed
to converge to a global minimum (since the residual is not globally convex),
and for this reason its limit may miss the optimaliy benchmark. On the other
hand, using the estimator derived in \S 3 as a ``good initialization
point'' to this mimimization algorithm leads to a limit state that has at least
the same order of accuracy.

These various approaches are numerically tested in \S 5 for the elliptic
equation \iref{ellip}, for both the overdetermined regime $m\geq d$, 
and the underdetermined regime $m<d$.

\section{Optimal recovery benchmarks}

In this section we describe the performance of
the best possible recovery map
\be
w\mapsto u^*(w),
\ee
in terms of its worst case error.  \dwo{We consider first the case
of noiseless data and no model bias.  In a subsequent step we   take}
such perturbations into account. While these
best recovery maps cannot be implemented by 
a simple algorithm, their performance serves as benchmark
for the nonlinear state estimation algorithms discussed in the
next section.

\subsection{Optimal recovery for the solution manifold}

In the absence of model bias and when a noiseless measurement $w=P_W u$ is given, our knowledge
on $u$ is that it belongs to the set
\be
\cM_w:=\cM\cap V_w.
\ee
The best possible recovery map can be described through the following general notion.

\begin{definition}
The Chebychev ball of a bounded set $S\in V$ is the closed ball $B(v,r)$ of minimal radius 
that contains $S$. One denotes by $v={\rm cen}(S)$ the Chebychev center of $S$ and 
$r={\rm rad}(S)$ its Chebychev radius. 
\end{definition}

\begin{figure}[ht]
\begin{subfigure}{.5\textwidth}
  \centering
  \includegraphics[scale=0.45]{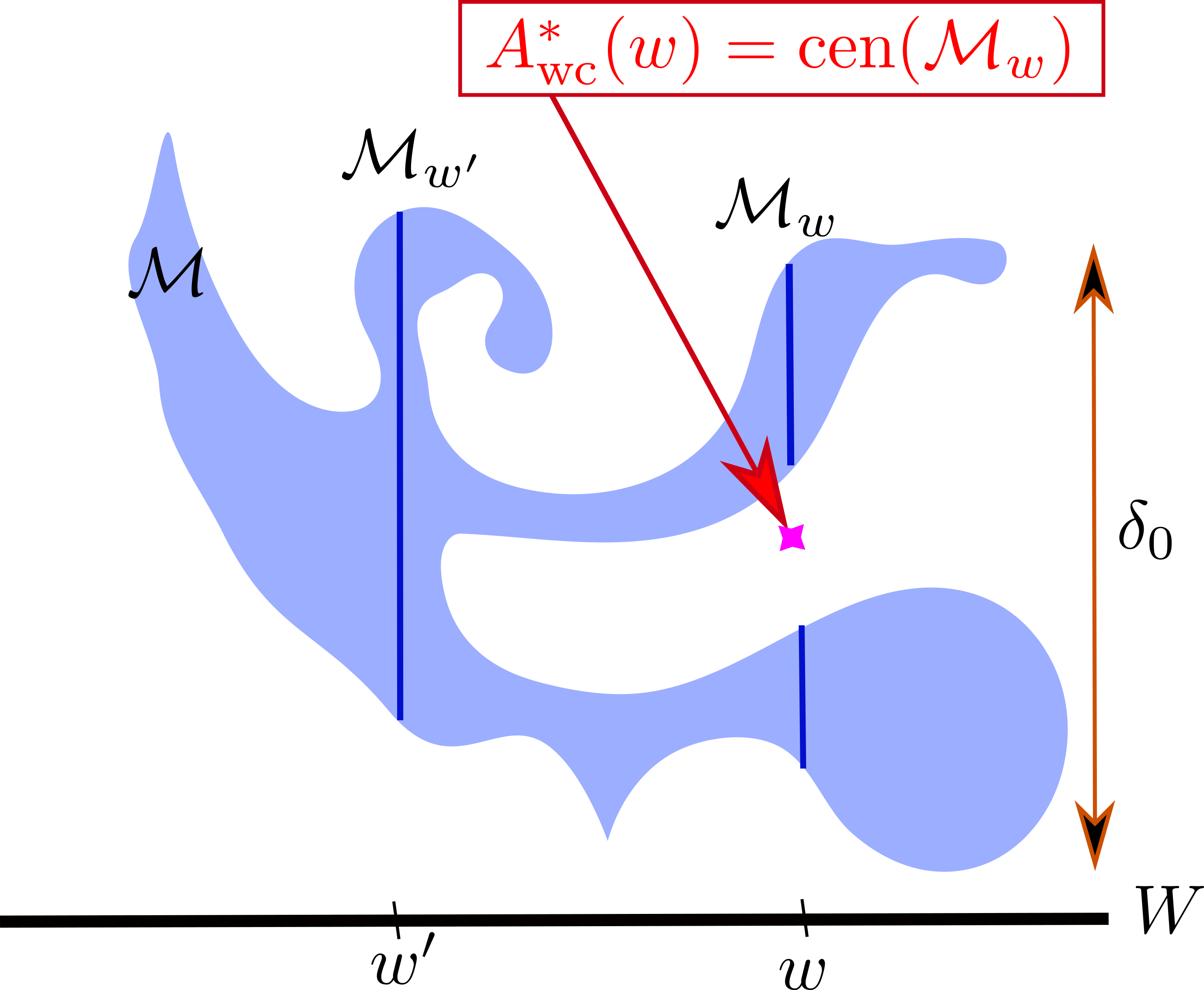}
  \caption{Perfect model.}
  \label{fig:bird}
\end{subfigure}
\begin{subfigure}{.5\textwidth}
  \centering
  \includegraphics[scale=0.45]{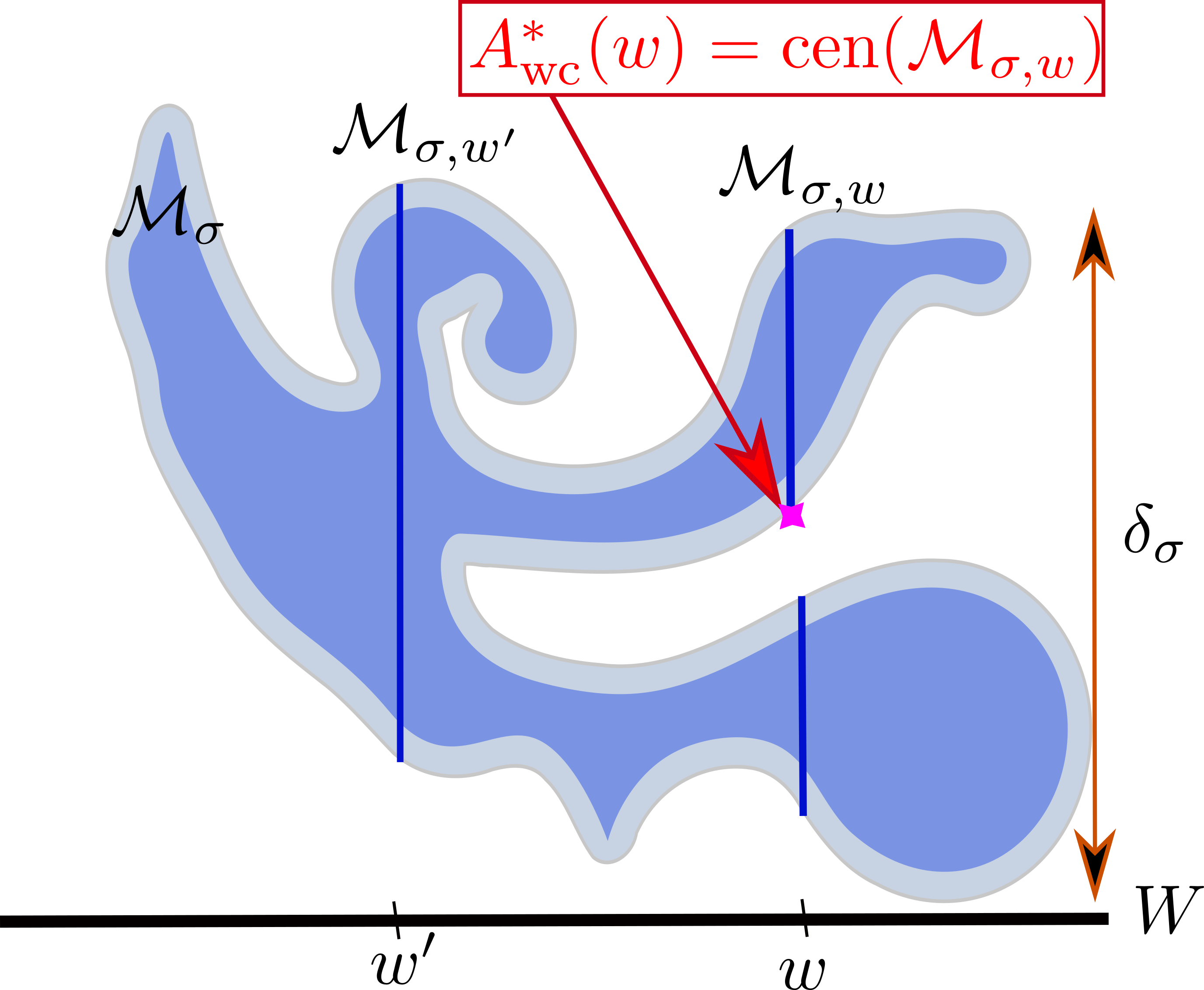}
  \caption{Model bias.}
  \label{fig:bird-offset}
\end{subfigure}
\caption{Illustration of the optimal recovery benchmark on a manifold in the two dimensional Euclidean space.}
\label{fig:fig}
\end{figure}

\ac{In particular 
one has
\be
\frac 1 2 {\rm diam}(S)\leq {\rm rad}(S)\leq  {\rm diam}(S),
\label{radiam}
\ee
where ${\rm diam}(S):=\ac{\sup}\{\|u-v\|\,:\, u,v\in S\}$ is the diameter of $S$.}
Therefore, the recovery map that minimizes the worst case error over $\cM_w$ for
any given $w$, and therefore over $\cM$ is
defined by
\be
u^*(w)={\rm cen}(\cM_w).
\ee
Its worst case error is
\be
E_\wc^*= \sup \{\rad(\cM_w)\, :\, w\in W\}.
\ee
In view of the equivalence \iref{radiam}, we can
relate $E^*_\wc$ to the quantity 
\be
\delta_0=\delta_0(\cM,W):=\ac{\sup}\{{\rm diam}(\cM_w)\,:\, w\in W\}=\sup \{\|u-v\|\; : \; u,v\in \cM, \;u-v\in W^\bot \},
\ee
by the equivalence
\be
\frac 1 2 \delta_0\leq E^*_\wc \leq \delta_0.
\ee
Note that injectivity of the measurement map $P_W$ over $\cM$
is equivalent to $\delta_0=0$. We provide in Figure \iref{fig:bird} an illustration the above benchmark concepts.

If $w=P_Wu$ for some $u\in \cM$, then any $u^*\in \cM$ such that $P_Wu^*=w$,
meets the ideal benchmark $\|u-u^*\|\leq \delta_0$. Therefore, one way
for finding such a $u^*$ would be to minimize the distance to the 
manifold over all functions such that $P_Wv=w$, that is, solve
\be
\min_{v\in V_w} {\rm dist} (v,\cM)=\min_{v\in V_w} \min_{y\in Y}\|u(y)-v\|.
\ee
This problem is computationally out of reach since it amounts
to  the nested minimization of two non-convex functions in high dimension.

Computationally feasible algorithms such as the PBDW methods
are based on a simplification of the manifold $\cM$ which induces an 
approximation error. We introduce \dwo{next a somewhat relaxed} benchmark
that takes this error into account.

\subsection{Optimal recovery under perturbations}

In order to account for manifold simplification as well as model bias,
for any given accucary $\sigma>0$, we introduce the $\sigma$-offset of $\cM$,
\be
\cM_\sigma:=\{v\in V\; : \; \dist(v,\cM) \leq \sigma\}=\bigcup_{u\in\cM} B(u,\sigma).
\ee
Likewise, we introduce the perturbed set 
\be
\cM_{\sigma,w} =\cM_\sigma\cap V_w,
\label{msigmaw}
\ee
\dwo{which, however, still excludes uncertainties in $w$.}
Our benchmark for the worst case error is now defined as (see Figure \iref{fig:bird-offset} for an illustration)
\be
\delta_\sigma:=\max_{w\in W} {\rm diam}(\cM_{\sigma,w})=\max \{\|u-v\|\; : \; u,v\in \cM_\sigma, \;u-v\in W^\bot \}.
\label{benchapp}
\ee
The map $\sigma\mapsto \delta_\sigma$ satisfies some elementary properties:
\begin{itemize} 
\item
Monotonicity and continuity: it is obviously non-decreasing
\be
\sigma\leq \t \sigma  \implies \delta_\sigma\leq \delta_{\t\sigma}.
\ee
Simple finite dimensional examples show that this map may have 
jump discontinuities. Take for example a compact set $\cM\subset \R^2$ consisting of the two 
points $(0,0)$ and $(1/2,1)$, and $W=\R e_1$ where $e_1=(1,0)$. Then
$\delta_\sigma=2\sigma$ for $0\leq \sigma \leq \frac 1 4$, while $\delta_{\frac 1 4}(\cM,W)=1$.
Using the compactness of $\cM$, it is possible to check that $\sigma\mapsto \delta_\sigma$
is continuous from the right and in particular $\lim_{\sigma \to 0} \delta_\sigma(\cM,W)=\delta_0$.
\item
Bounds from below and above: for any \om{$u,v\in \cM_{\sigma, w}$}, and for any $\t \sigma\geq 0$,
\dwo{let $\t u=u+\t \sigma g$ and $\t v=v-\t \sigma g$ with $g=(u-v)/\|u-v\|$.
Then,}  $\|\t u-\t v\|=\|u-v\|+2\om{\tilde\sigma}$ and $\t u-\t v\in W^\bot$, \dwo{which shows that} \om{$\t u,\, \t v\in \cM_{\sigma+\t \sigma, w}$,  and}
\be
\delta_{\sigma+\t \sigma} \geq \delta_{\sigma} +2\t \sigma.
\ee 
In particular,
\be
\delta_\sigma \geq \delta_0 + 2\sigma\geq 2\sigma.
\label{2sigma}
\ee
On the other hand, we obviously have the upper bound $\delta_\sigma\leq {\rm diam}(\cM_\sigma)\leq {\rm diam}(\cM)+2\sigma$.
\item
The quantity
\be
\mu(\cM,W):=\frac 1 2\sup_{\sigma> 0}\frac {\delta_\sigma-\delta_0}{\sigma},
\label{globstab}
\ee
may be viewed as a general stability constant inherent to the recovery problem, similar to $\mu(V_n,W)$
that is more specific to the particular PBDW method: in the special case where $\cM=V_n$
and $V_n\cap W^\perp=\{0\}$, one has $\delta_0=0$ and $\frac {\delta_\sigma}{\sigma}=\mu(V_n,W)$ for all $\sigma>0$.
Note that $\mu(\cM,W)\geq 1$ in view of \iref{2sigma}.
\end{itemize}

Regarding measurement noise, it suggests to introduce the quantity
\be
\t \delta_\sigma:=\max \{\|u-v\|\; : \; u,v\in \cM,\; \|P_W u-P_W v\|\leq \sigma \}.
\ee
\ac{The two quantities $\delta_\sigma$ and $\t \delta_\sigma$ are not equivalent, however}
one has the framing
\be
\delta_\sigma-2\sigma\leq \t \delta_{2\sigma}\leq \delta_\sigma+2\sigma.
\label{framing}
\ee
To prove the upper inequality, we note that for any $u,v\in \cM_\sigma$ such that $u-v\in W^\bot$, there
exists $\t u,\t v\in\cM$ at distance $\sigma$ from $u,v$, respectively, and therefore such that
$\|P_W(\t u\dw{-} \t v)\| \leq 2\sigma$. Conversely, for any $\t u,\t v\in\cM$ such that $\|P_W(\t u\dw{-} \t v)\| \leq 2\sigma$, there exists
$u,v$ at distance $\sigma$ from $\t u,\t v$, respectively, such that $u-v\in W^\bot$,
which gives the lower inequality.
Note that the upper inequality in \iref{framing} combined with \iref{2sigma}
implies that $\t \delta_{2\sigma}\leq 2 \delta_{\sigma}$.

In \dwo{the following}  analysis of reconstruction methods,
we use the quantity $\delta_\sigma$ as a benchmark which, in view of this last observation, also
accounts for the lack of accuracy in the measurement of $P_Wu$.
Our objective is therefore to design an algorithm that, for a given tolerance $\sigma>0$,
recovers from the measurement $w=P_Wu$ an approximation to $u$
with accuracy comparable to $\delta_\sigma$. Such an algorithm requires 
that we are able to capture the solution manifold up to some tolerance 
$\e\leq \sigma$ by some reduced model. 

\section{Nonlinear recovery by reduced model selection}

\subsection{Piecewise affine reduced models}

Linear or affine reduced models, as used in the PBDW algorithm,
are not suitable for approximating the solution manifold
when the required tolerance $\e$ is too small.
In particular, when $\e\dwo{<} d_m(\cM)$ one would then need to use a linear space $V_n$ of 
dimension $n\dwo{>} m$, therefore \ac{making} $\mu(V_n,W)$ infinite.

One way out is to replace the single space $V_n$ by a 
{\em family} of affine spaces 
\be
V_k=\o u_k+\o V_k, \quad k=1,\dots,K,
\ee
each of them having dimension 
\be
\dim(V_k)=n_k\leq m,
\ee
such that the manifold is well captured by the union of these spaces,
in the sense that 
\be
{\rm dist}\(\cM,\bigcup_{k=1}^K V_k\)\leq \e
\ee
for some prescribed tolerance $\e>0$. This is equivalent to saying that there exists a partition
of the solution manifold
\be
\cM=\bigcup_{k=1}^K \cM_k,
\ee
such that we have local certified bounds
\be
{\rm dist}(\cM_k,V_k)\leq \e_k \leq \e,\quad k=1,\dots,K.
\label{localacc}
\ee
We may thus think of the family $(V_k)_{k=1,\dots,K}$ as a piecewise
affine approximation to $\cM$. We stress that, in contrast to the hierarchies 
$(V_n)_{n=0,\dots,m}$ of reduced models discussed in \S 1.2, the spaces $V_k$
do not have dimension $k$ and are not nested. \dwo{Most importantly,} $K$ is not limited
by $m$ while each $n_k$ is.

The objective of using a piecewise reduced model in the context of state
estimation is to have a joint control on the local accuracy $\e_k$
as expressed by \iref{localacc} and on the stability of the PBDW when using any
individual $V_k$. This means that, for some prescribed 
$\mu>1$, we ask that
\be
\mu_k=\mu(\o V_k,W) \leq \mu, \quad k=1,\dots,K.
\label{localstab}
\ee
According to \iref{rec_error}, the worst case error bound over $\cM_k$ when using the PBDW method
with a space $V_k$ is given by the product $\mu_k\e_k$. This suggests to 
alternatively require from the collection $(V_k)_{k=1,\dots,K}$, that for some 
prescribed $\sigma>0$, one has
\be
\sigma_k:=\mu_k\e_k \leq \sigma, \om{ \quad k=1,\dots,K.}
\label{sigmaadm}
\ee
\ac{
This leads us to the following definitions.}

\begin{definition}
\ac{The family $(V_k)_{k=1,\dots,K}$ is {\em $\sigma$-admissible}
if \iref{sigmaadm} holds. It is {\em $(\e,\mu)$-admissible}
if \iref{localacc} and \iref{localstab} are jointly satisfied.}
\end{definition}

Obviously, any $(\e,\mu)$-admissible
family is $\sigma$-admissible with $\sigma:=\mu\e$.
In this sense the notion of $(\e,\mu)$-admissibility is thus more restrictive
than that of $\sigma$-admissibility. The benefit of the first notion is in 
the uniform control on the size of $\mu$ which is
critical in the presence of noise, as hinted at by
Remark \ref{rembiasnoise}.

If $u\in \cM$ is our unknown state and $w=P_W u$ is its observation, we may apply the 
PBDW method for the different $V_k$ in the given family, which yields a corresponding 
family \dw{of} estimators
\be
u_k^*=u_k^*(w)={\rm argmin}\{\dist(v,V_k)\, : \, v\in V_w\}, \quad k=1,\dots,K.
\label{estimators}
\ee
If $(V_k)_{k=1,\dots,K}$ is $\sigma$-admissible, we find that the accuracy bound
\be
\|u-u_k^*\|\leq \mu_k{\rm dist}(u,V_k) \leq \mu_k\e_k= \sigma_k \leq \sigma,
\ee
holds whenever $u\in \cM_k$.

Therefore, if in addition to the observed data $w$
one had an oracle giving the information on which portion $\cM_k$ of the manifold the unknown state sits,
we could derive an estimator with worst case error 
\be
E_{wc}\leq \sigma. 
\label{oracle}
\ee
This information is, however, not available and such a worst case error estimate
cannot be hoped for, even with an additional multiplicative constant. Indeed, 
as we shall see below, $\sigma$ can be fixed arbitrarily small
by the user when building the family $(V_k)_{k=1,\dots,K}$, while we know from \S 2.1 that the
worst case error is bounded from below by $E_{wc}^*\geq \frac 1 2 \delta_0$ which could be non-zero.
We will thus need to replace the ideal choice of $k$ by a model selection procedure
only based on the data $w$, that is, a map
\be
w \mapsto k^*(w),
\ee
leading to a choice of estimator $u^*=u^*_{k^*}$. We shall prove further that such an estimator
is able to achieve the accuracy
\be
E_{wc} \leq \delta_{\sigma},
\ee
that is, the benchmark introduced in \S 2.2. Before discussing this model selection, we discuss the 
existence and construction of $\sigma$-admissible or $(\e,\mu)$-admissible families.

\subsection{Constructing admissible reduced model families} \label{ssec:admissible_families}

For any arbitrary choice of $\e>0$ and $\mu\geq 1$, the existence of an $(\e,\mu)$-admissible family
results from the following observation: since the manifold $\cM$ is a compact set of $V$, there exists
a finite $\e$-cover of $\cM$, that is, a family $\o u_1,\dots,\o u_K \in V$ such that
\be
\cM\subset \bigcup_{k=1}^K B(\o u_k,\e),
\ee
or equivalently, for all $v\in \cM$, there exists a $k$ such that $\|v-u_k\|\leq \e$. With such an $\e$ cover,
we consider the family of trivial affine spaces defined by
\be
V_k=\{\o u_k\}=\o u_k+\o V_k, \quad \o V_k=\{0\},
\ee
thus with $n_k=0$ for all $k$. The covering property implies that \iref{localacc} holds. On the other hand,
for the $0$ dimensional space, one has
\be
\mu(\{0\},W)=1,
\ee
and therefore \iref{localstab} also holds. The family $(V_k)_{k=1,\dots,K}$ is therefore $(\e,\mu)$-admissible,
and also  $\sigma$-admissible with $\sigma=\e$.

This family is however not satisfactory for algorithmic purposes for two main reasons. First, the manifold
is not explicely given to us and the construction of the centers $\o u_k$ is by no means trivial. Second, 
asking for an $\e$-cover, would typically require that $K$ becomes extremely large as $\e$ goes to $0$. 
For example, assuming that the parameter to solution $y\mapsto u(y)$ has Lipschitz constant $L$,
\be
\|u(y)-u(\t y)\|\leq L|y-\t y|, \quad y,\t y\in Y,
\ee
for some norm $|\cdot|$ of $\R^d$, then an $\e$ cover for $\cM$ would be induced by an $L^{-1}\e$ cover for $Y$
which has cardinality $K$ growing like $\e^{-d}$ as $\e\to 0$. Having a family of moderate size $K$ is important
for the estimation procedure since we intend to apply the PBDW method for all $k=1,\dots,K$. 

In order to construct $(\e,\mu)$-admissible or $\sigma$-admissible families of better controlled size, we need to split the
manifold in a more economical manner than through an $\e$-cover, and use spaces $V_k$ of general
dimensions $n_k\in \{0,\dots,m\}$ for the various manifold portions $\cM_k$. To this end, we combine
standard constructions of linear reduced model spaces with an iterative splitting
procedure operating on the parameter domain $Y$. Let us mention that various ways 
of splitting the parameter domain have already been considered in order to produce local 
reduced bases having both, controlled cardinality and prescribed accuracy \cite{EPR,MS,BCDGJP}. Here our goal
is slightly different since we want to control both the accuracy $\e$ and the stability $\mu$ with respect to 
the measurement space $W$.

We describe the greedy algorithm for constructing $\sigma$-admissible families, and explain how it
should be modified for $(\e,\mu)$-admissible families. For simplicity \dw{we} consider the case where
$Y$ is a rectangular domain with sides parallel to the main axes, the extension to a more general
bounded domain $Y$ being done by embedding it in such a \dw{hyper}-rectangle. We are given a prescribed target value $\sigma>0$
and the splitting procedure starts from $Y$. 

At step $j$,
a disjoint partition of $Y$ into rectangles $(Y_k)_{k=1,\dots,K_j}$ with sides parallel
to the main axes has been generated. It induces
a partition of $\cM$ given by
\be
\cM_k:=\{u(y)\,:\, y\in Y_k\}, \quad k=1,\dots,K_j.
\ee
To each $k\in \{1,\dots,K_j\}$ we associate a hierarchy of affine reduced basis spaces
\be \label{eq:aff_from_basis}
V_{n,k}=\o u_k+\o V_{n,k}, \quad n=0,\dots, m.
\ee
where $\o u_k=u(\o y_k)$ with $\o y_k$ the vector defined as the center of the rectangle $Y_k$. The nested linear 
spaces
\be
\o V_{0,k}\subset \o V_{1,k}\subset \cdots\subset \o V_{m,k}, \quad \dim(\dwo{\o V}_{n,k})=n,
\ee
are meant to approximate the translated portion of the manifold $\cM_k-\o u_k$. 
For example, they could be reduced basis spaces obtained by applying the 
greedy algorithm to $\cM_k-\o u_k$, or spaces resulting from local $n$-term
polynomial approximations of $u(y)$ on the rectangle  $Y_k$.
Each space $V_{n,k}$ has a given accuracy bound
and stability constant 
\be \label{eq:acc_and_stab}
{\rm dist}(\cM_k,V_{n,k}) \leq  \e_{n,k}  \quad{\rm and}\quad \mu_{n,k}:=\mu(\o V_{n,k},W).
\ee 

We define the test quantity
\be \label{eq:pm_criteria}
\tau_k=\min_{n=0,\dots,m} \mu_{n,k}\e_{n,k}.
\ee
If $\tau_k\leq \sigma$, the rectangle $Y_k$ is not \dwo{split} and becomes a member of the final partition.
The affine space associated to $\cM_k$ is 
\be
V_k=\o u_k+\o V_k,
\ee
where $V_k=V_{n,k}$ for the value of $n$ that minimizes $\mu_{n,k}\e_{n,k}$. The rectangles
$Y_k$ with $\tau_k> \sigma$ are, on the other hand, \dwo{split} into a finite number of sub-rectangles 
in a way that we discuss below. This results in the new larger partition 
$(Y_k)_{k=1,\dots,K_{j+1}}$ after relabelling the $Y_k$. The algorithm terminates
at the step $j$ as soon as $\tau_k\leq \sigma$ for all $k=1,\dots,K_j=K$,
and the family $(V_k)_{k=1,\dots,K}$ is $\sigma$-admissible.
In order to obtain an $(\e,\mu)$-admissible family, we simply modify the test quantity $\tau_k$
by defining it instead as
\be
\tau_k:=\min_{n=0,\dots,m}\max\Big\{\frac {\mu_{n,k}}\mu, \frac {\e_{n,k}}\e\Big\}
\ee
and splitting the cells for which $\tau_k>1$.

The splitting of one single rectangle $Y_k$ can be performed in various ways. When the parameter 
dimension $d$ is moderate, we may subdivide each side-length at the mid-point, resulting into $2^d$ sub-rectangles
of equal size. This splitting becomes too   costly as $d$ gets large, in which case it is preferable to make a choice 
of $i\in \{1,\dots,d\}$ and subdivide $Y_k$ at the mid-point of the side-length in the $i$-coordinate, resulting
in  only $2$ sub-rectangles. In order to decide which coordinate to pick, we consider the $d$ possibilities and
take the value of $i$ that minimizes the quantity
\be
\tau_{k,i}=\max \{\tau_{k,i}^{\dwo{-}},\tau_{k,i}^+\},
\ee
where $(\tau_{k,i}^{\dwo{-}},\tau_{k,i}^+)$ are the values of $\tau_k$ for the two subrectangles obtained by splitting along the $i$-coordinate.
In other words, we split in the direction that decreases $\tau_k$  most effectively. In order to be certain that all sidelength
are eventually \dwo{split}, we can mitigate the greedy choice of $i$ in the following way: if $Y_k$ has been generated
by $l$ consecutive refinements, and therefore has volume $|Y_k|=2^{-l}|Y|$, and if $l$ is even, we choose $i=\ac{(l/2 \,{\rm mod}\,d)}$. 
This 
means that at each even level we split in a cyclic manner in the coordinates $i\in\{1,\dots,d\}$.

Using such elementary splitting rules, we are ensured that the algorithm must terminate. Indeed, we are guaranteed
that for any $\eta>0$, there exists a level $l=l(\eta)$ such that any rectangle $Y_k$ generated by $l$ consecutive refinements
has side-length smaller than $2\eta$ in each direction. Since the parameter-to-solution map is continuous, for any $\e>0$, we can pick
$\eta>0$ such that
\be
\|y-\t y\|_{\ell^\infty}\leq \eta \implies \|u(y)-u(\t y)\|\leq \e, \quad y,\t y\in Y.
\ee
Applying this to $y\in Y_k$ and $\t y=\o y_k$, we find that \dwo{for $\o u_k = u(\o y_k)$}
\be
\|u-\o u_k\| \leq \e, \quad u\in \cM_k.
\ee
Therefore, for any rectangle $Y_k$ of generation $l$, we find that the trivial affine space $V_k=\o u_k$ has local accuracy $\e_k\leq \e$
and $\mu_k=\mu(\{0\},W)=1\leq \mu$, which implies that such a rectangle would not anymore be refined by the algorithm.

\subsection{Reduced model selection and recovery bounds} \label{ssec:recovery_bounds}
\label{sec:model-selection}
We return to the problem of selecting an estimator 
within the family $(u_k^*)_{k=1,\dots,K}$ defined by \iref{estimators}.
In an idealized version, the selection procedure picks the value $k^*$
that minimizes the distance of $u_k^*$ to the solution manifold, that is,
\be
k^*={\rm argmin}\{ \dist(u_k^*,\cM)\,:\, k=1,\dots,K\}
\label{idealsel}
\ee
and takes for the final estimator
\be
\label{ustar}
u^*=u^*(w):=u^*_{k^*}(w).
\ee
Note that $k^*$ also depends on the observed data $w$.
This estimation procedure is not realistic, since the computation of the distance
of a known function $v$ to the manifold
\be
\dist(v,\cM)=\min_{y\in Y} \|u(y)-v\|,
\ee
is a high-dimensional non-convex  problem, which
necessitates to explore the solution manifold. A more realistic 
procedure is based on replacing this distance by a
surrogate quantity 
${\cal S}(v,\cM)$ that is easily computable
and satisfies a uniform equivalence
\be
r\dist(v,\cM)\leq {\cal S}(v,\cM)\leq R \dist(v,\cM), \quad v\in V,
\label{surframe}
\ee
for some constants $0<r\leq R$. We then instead take for $k^*$ the value
that minimizes this surrogate, that is,
\be
k^*={\rm argmin}\{ {\cal S}(u_k^*,\cM)\,:\, k=1,\dots,K\}.
\label{realisticsel}
\ee
Before discussing the derivation of ${\cal S}(v,\cM)$ in concrete cases, we establish a recovery
bound in the absence of model bias and noise.

\begin{theorem}
\label{theorecovery}
Assume that the family $(V_k)_{k=1,\dots,K}$ is $\sigma$-admissible for some $\sigma>0$. Then, the idealized estimator
based on \iref{idealsel}, \eref{ustar}, satisfies the worst case error estimate
\be
E_{wc}=\max_{u\in \cM} \|u-u^*(P_Wu)\| \leq \delta_\sigma,
\label{recovideal}
\ee
where $\delta_\sigma$ is the benchmark quantity defined in \iref{benchapp}. When using the estimator
based on \iref{realisticsel}, the worst case error estimate is modified into
\be
E_{wc} \leq \delta_{\kappa\sigma}, \quad \kappa=\frac R r>1.
\label{recovreal}
\ee
\end{theorem}

\noindent
{\bf Proof:} Let $u\in \cM$ be an unknown state and $w=P_Wu$.
There exists $l=l(u)\in {1,\dots,K}$, such that $u\in\cM_{l}$,
and for this value, we know that
\be
\|u-u_{l}^* \| \leq \mu_{l}\e_{l}=\sigma_{l}\leq \sigma.
\ee
Since $u\in \cM$, it follows that
\be
\dist(u_{l}^*,\cM)\leq \sigma.
\ee
On the other hand, for the value $k^*$ selected by \iref{realisticsel} and $u^*=u_{k^*}^*$, we have
\be
\dist(u^*,\cM)\leq R\, {\cal S}(u^*,\cM)\leq R\, {\cal S}(u_{l}^*,\cM) \leq \kappa \dist(u_{l}^*,\cM) \leq \kappa\sigma.
\ee
It follows that $u^*$ belongs to the offset $\cM_{\om{\kappa}\sigma}$. Since $u\in \cM\subset \cM_\sigma\om{\subseteq \cM_{\kappa\sigma}}$ and $u-u^*\in W^\perp$,
we find that 
\be
\|u-u^*\|\leq \delta_{\kappa\sigma},
\ee
which establishes the recovery estimate \iref{recovreal}. The estimate  \iref{recovideal} for the idealized estimator follows since it 
corresponds to having $r=R=1$. \hfill $\Box$

\begin{remark}
One possible variant of the selection mechanism, which is actually adopted in our numerical
experiments, consists in picking the value $k^*$
that minimizes the distance of $u_k^*$ to the corresponding local portion $\cM_k$
of the solution manifold, or a surrogate $\cS(u_k^*,\cM_k)$ with 
equivalence properties analogous to \iref{surframe}. It is readily checked that Theorem \ref{theorecovery} remains valid for the resulting estimator $u^*$ with the same type of proof.
\end{remark}

In the above result, we do not obtain the best possible accuracy 
satisfied by the different $u_k^*$, since we do not have
an oracle providing the information on the best choice of $k$. We next show that this order of 
accuracy is attained in the particular case where the measurement map 
$P_W$ is injective on $\cM$ and the stability constant of the recovery
problem defined in \iref{globstab} is finite.

\begin{theorem}
\label{thm:3.3}
Assume that $\delta_0=0$ and that
\be
\mu(\cM,W)=\frac 1 2\sup_{\sigma>0} \frac{\delta_\sigma}{\sigma}<\infty.
\ee
Then, for any given state $u\in \cM$ with observation $w=P_Wu$,
the estimator $u^*$ obtained by the model selection procedure \iref{realisticsel} satisfies
the oracle bound
\be
\|u-u^*\| \leq C\min_{k=1,\dots,K} \|u-u_k^*\|, \quad C:=2\mu(\cM,W)\kappa.
\label{oracle1}
\ee 
In particular, if $(V_k)_{k=1,\dots,K}$ is $\sigma$-admissible, it satisfies
\be
\|u-u^*\| \leq C\sigma.
\label{oracle2}
\ee
\end{theorem}

\noindent
{\bf Proof:} Let $l\in \{1,\dots,K\}$ be the value for which $\|u-u_l^*\|=\min_{k=1,\dots,K} \|u-u_k^*\|$. 
Reasoning as in the proof of Theorem \ref{theorecovery}, we find that 
\be
\dist(u^*,\cM)\leq \kappa \beta, \quad \beta:= \dist(u_{l}^*,\cM),
\ee
and therefore
\be
\|u-u^*\|\leq \delta_{\kappa\beta}\leq 2 \mu(\cM,W) \kappa \dist(u_{l}^*,\cM),
\ee
which is \iref{oracle1}. We then obtain \iref{oracle2} using the fact that 
$\|u-u^*_k\|\leq \sigma$ for the value of $k$ such that $u\in \cM_k$.
\hfill $\Box$
\newline

We next discuss how to incorporate model bias and noise in the recovery bound, provided that we have
a control on the stability of the PBDW method, through a uniform bound on $\mu_k$, which holds
when we use $(\e,\mu)$-admissible families.

\begin{theorem}
\label{thm:3.4}
Assume that the family $(V_k)_{k=1,\dots,K}$ is $(\e,\mu)$-admissible for some $\e>0$ and $\mu\geq 1$. 
If the observation is $w=P_Wu+\eta$ with $\|\eta\|\leq \e_{noise}$,
and if the true state does not lie in $\cM$ but satisfies ${\rm dist}(u,\cM)\leq \e_{model}$,
then, the estimator based on \iref{realisticsel} satisfies the estimate
\be
\|u-u^*(w)\|\leq \delta_{\kappa\rho}+\e_{noise}, \quad \rho:=\mu(\e+\e_{noise})+(\mu+1)\e_{model},\quad \kappa=\frac R r,
\label{pertest}
\ee
and the idealized estimator based on \iref{idealsel} satifies a similar estimate with $\kappa=1$.
\end{theorem}

\noindent
{\bf Proof:} There exists $l=l(u)\in \{1,\dots,K\}$ such that 
\be
\dist(u,\cM_{l})\leq \e_{model},
\ee 
and therefore
\be
\dist(u,V_{l})\leq \e_{l}+\e_{model}.
\ee 
As already noted in Remark \ref{rembiasnoise}, we know that the PBDW method for this value of $k$ has accuracy
\be
\|u-u_{l}^* \| \leq \mu_{l}(\e_{l}+\e_{noise}+\e_{model})\leq \mu(\e+\e_{noise}+\e_{model}).
\ee
Therefore
\be
\dist(u_{l}^*,\cM)\leq  \mu(\e+\e_{noise}+\e_{model})+\e_{model}=\rho,
\ee
and in turn
\be
\dist(u^*,\cM)\leq \kappa \rho.
\ee
On the other hand, we define 
\be
v:=u+\eta=u+w-P_W u=u+P_W(u^*-u), 
\ee
so that 
\be
\dist(v,\cM)\leq \|v-u\|+\e_{model} \leq \e_{noise}+\e_{model} \leq \rho.
\ee
Since $v-u^*\in W^\perp$, we conclude that $\|v-u^*\| \leq \delta_{\kappa\rho}$, from which \iref{pertest} follows.
\hfill $\Box$
\nl

While the reduced model selection approach provides \dwo{us with an estimator $w\mapsto u^*(w)$ of a single 
plausible state, the estimated distance of some of the other estimates $u_k(w)$ may be of comparable size. Therefore,}
one could be interested in recovering a more complete estimate on \dwo{a} plausible {\em set} \dwo{that may contain   the true
state $u$ or even several states in $\cM$ sharing the same measurement.} This more ambitious goal can be viewed as a deterministic counterpart to the search for the 
entire posterior probability distribution of the state in a Bayesian estimation framework, instead of
only searching for \dwo{a single  estimated state, for instance,} the expectation of this distribution. For simplicity, 
we discuss this problem in the absence of model bias and noise. Our goal is therefore to approximate
the set 
\be
\cM_w=\cM \cap V_w.
\ee
Given the family $(V_k)_{k=1,\dots,K}$, we consider the ellipsoids
\be
\cE_k:=\{v\in V_w \,\; \, \dist(v,V_k)\leq \e_k\}, \quad k=1,\dots,K,
\ee 
which have center $u^*_k$ and diameter \dwo{at most} $\mu_k\e_k$. We already 
know that $\cM_w$ is contained inside the union of the $\cE_k$ which could
serve as a first estimator. \ac{In order to refine this estimator, we would like
to discard the $\cE_k$ that do not intersect the
associated portion $\cM_k$ of \dw{the} solution manifold.}

For this purpose, we define our estimator of $\cM_w$ as the union
\be
\cM_w^*:=\bigcup_{k\in S} \cE_k,
\ee
where $S$ is the set of those $k$ such that
\be
{\cal S}(u_k^*,\cM_{\dwo{k}}) \leq R\mu_k\e_k.
\ee
It is readily seen that $k\notin S$ implies that $\cE_k\cap \cM_k = \emptyset$.
The following result shows that this set approximates $\cM_w$ with an accuracy of the same order
as the recovery bound established for the estimator \dwo{$u^*(w)$}.

\begin{theorem}
\label{theorecoveryset}
For any state $u\in \cM$ with observation $w=P_Wu$, one has the inclusion
\be
\cM_w \subset \cM_w^*.
\label{incl}
\ee
If the family $(V_k)_{k=1,\dots,K}$ is $\sigma$-admissible for some $\sigma>0$, 
the Hausdorff distance between the two sets satisfies the bound
\be
d_{\cH}(\cM_{w}^*,\cM_w)=\max_{v\in \cM_{w}^*}\min_{u\in \cM_w} \|v-u\|\leq \delta_{(\kappa+1)\sigma},\quad \kappa=\frac R r.
\label{setbound}
\ee
\end{theorem}

\noindent
{\bf Proof:} Any $u\in \cM_w$ is a state from $\cM$ that gives the observation $P_Wu$.
This state belongs to $\cM_l$ for some particular $l=l(u)$,
for which we know that $u$ belongs to the ellipsoid $\cE_l$
and that
\be
\| u-u^*_l\| \leq \mu_l\e_l.
\ee
This implies that $\dist(u^*_l,\cM_{\dw{l}})\leq \mu_l\e_l$, and therefore ${\cal S} (u^*_l,\cM_{\dw{l}})\leq R\mu_l\e_l$.
Hence $l\in S$, which proves the inclusion \iref{incl}. In order to prove
the estimate on the Hausdorff distance, we take any $k\in S$, and notice that
\be
\label{distk}
\dist(u_k^*,\cM_{\dw{k}})\leq \kappa \mu_k\e_k \leq \kappa\sigma,
\ee
and therefore, for all such $k$ and all $v\in \cE_k$, we have
\be
\dist(v,\cM_{\dw{k}})\leq (\kappa+1)\mu_k\e_k.
\ee
Since $u-v\in W^\perp$, it follows that
\be
\|v-u\|\leq  \delta_{(\kappa+1)\sigma},
\ee
which proves \iref{setbound}.
\hfill $\Box$

\begin{remark}
\ac{If we could take $S$ to be exactly the set of those $k$
such that $\cE_k\cap \cM_k \neq \emptyset$, the resulting $\cM^*_w$
would still contain $\cM_w$ but with a sharper error bound.
Indeed, any $v\in \cM^*_w$ belongs
to a set $\cE_k$ that intersects $\cM_k$ at some $u\in\cM_w$, so that
\be
d_{\cH}(\cM_{w}^*,\cM_w)\leq 2\sigma.
\ee
In order to identify if a $k$ belongs to this smaller $S$, we need to solve the
minimization problem
\be
\min_{v\in \cE_k}\cS(v,\cM_k),
\ee
and check whether the minimum is zero. As  explained next, the quantity 
$\cS(v,\cM_k)$ is itself obtained by a minimization problem over $y\in Y_k$.
The resulting double minimization problem is globally non-convex, but it is
convex separately in $v$ and $y$, which allows one to apply simple alternating minimization
techniques. These procedures (which are not guaranteed to
converge to the global minimum) are discussed in \S \ref{ssec:alternate}.}
\end{remark}

\subsection{Residual based surrogates}\label{ssec:residual}
\newcommand{\rZ}{r_{\mbox{\tiny Z}}}
\newcommand{\rV}{r_{\mbox{\tiny V}}}
\dwo{The computational  realization of the above concepts hinges on two main constituents, namely (i)  the ability to evaluate bounds $\e_n$ for 
$\dist (\cM,V_n)$ as well as  (ii)  to have at hand computationally affordable surrogates ${\cal S}(v,\cM)$ for $\dist (v,\cM)= \min_{u\in \cM}\|v -u\|$.
In both cases one exploits the fact that errors in $V$ are equivalent to residuals in a suitable dual norm. Regarding (i), the derivation of bounds
$\e_n$ has been discussed extensively in the context of Reduced Basis Methods \cite{RHP}, see also \cite{DPW} for
the more general framework discussed below. Substantial   computational effort in an offline phase 
provides residual based surrogates for  $\|u-u(y)\|$  permitting frequent parameter queries at an online stage needed, in particular, to construct
reduced bases. This strategy becomes challenging though for high parameter dimensionality and we refer to \cite{CDDN} for remedies based
on trading deterministic certificates against probabilistic ones at significantly reduced computational cost. Therefore, we focus here on task (ii).
}

One typical setting where a computable surrogate $\cS(v,\cM)$ can be derived is when
$u(y)$ is the solution to a parametrized operator equation of the general form
\be
A(y)u(y)=f(y).
\ee
Here we assume that for every $y\in Y$ the right side $f(y)$ belongs
to the dual $Z'$ of a Hilbert test space $Z$, and $A(y)$ is boundedly invertible from $V$ to $Z'$.
The operator equation has therefore an equivalent variational formulation
\be
\label{varform}
\cA_y(u(y),v)=\cF_y(v), \quad v\in Z,
\ee
with parametrized bilinear form $\cA_y(w,v)=\<A(y)w,v\>_{Z',Z}$ and linear form $\cF_y(v)=\<f(y),v\>_{Z',Z}$.
This setting includes classical elliptic problems with $Z=V$,
as well as saddle-point and \dwo{unsymmetric problems such as convection-diffusion problems or space-time formulations of parabolic problems}.

We assume continuous dependence of $A(y)$ and $f(y)$ with respect to $y\in Y$, which by compactness of $Y$,
implies uniform boundedness and invertibility, that is
\be
\|\cA(y)\|_{V\to Z'}\leq R\quad {\rm and}\quad \|\cA(y)^{-1}\|_{Z'\to V}\leq r^{-1}, \quad y\in Y,
\label{boundA}
\ee
for some $0<r\leq R<\infty$. It follows that for any $v\in V$, one has the equivalence
\be
\label{rR}
r\|v-u(y)\|_V \leq \cR(v,y) \leq R\|v-u(y)\|_V.
\ee
where
\be
\cR(v,y):=\|A(y)v-f(y)\|_{Z'},
\label{resid}
\ee 
is the residual of the PDE for a state $v$ and parameter $y$.

Therefore the quantity
\be
\label{surrS}
\cS(v,\cM):=\min_{y\in Y}\cR(v,y),
\ee
provides \dwo{us} with a surrogate of $\dist(v,\cM)$ that satisfies the required framing \iref{surframe}.

One first advantage of this surrogate quantity is that for each given $y\in Y$,
the evaluation of the residual $\|A(y)v-f(y)\|_{Z'}$ does not require to compute the solution $u(y)$.
Its second advantage is that the minimization in $y$ is facilitated in the relevant case where
$A(y)$ and $f(y)$ have affine dependence on $y$, that is,
\be
A(y)=A_0 +\sum_{j=1}^d y_jA_j\quad {\rm and} \quad f(y)=f_0+\sum_{j=1}^d y_j f_j.
\label{affine}
\ee
Indeed, $\cS(v, \cM)$ then amounts to the minimization over $y\in Y$ of the function 
\be
\label{Rvy0}
\cR(v,y)^2:=\Big\|A_0v-f_0+\sum_{j=1}^d y_j(A_j v-f_j)\Big\|_{Z'}^2,
\ee
which is a convex quadratic polynomial in $y$. \dwo{Hence, a minimizer $y(v)\in Y$ of the corresponding constrained linear least-squares problem exists, rendering the surrogate $\cS(v,\cM)= \cR(v,y(v))$ well-defined.}

\ac{In all the above mentioned examples the norm $\|\cdot\|_Z = \langle\cdot,\cdot\rangle_Z^{1/2}$ can be efficiently computed.  For instance, in the simplest case of an $H^1_0(\Omega)$-elliptic problem 
one has $Z=V = H^1_0(\Omega)$ with 
\be
\langle v,z\rangle_Z = \int_\Omega \nabla v\cdot\nabla z dx.
\ee 
The obvious obstacle is then, however, the computation of the {\em dual norm} $\|\cdot\|_{Z'}$ which in the particular example  above is the $H^{-1}(\Omega)$-norm. A viable strategy is to use} \ac{the Riesz lift $\rZ: Z'\to Z$, defined by
\be
\label{Riesz}
\langle \rZ g,z\rangle_Z = \<g,z\>_{Z',Z}= g(z),\quad g\in Z',\, z\in Z,
\ee
\dw{which implies} that $\|\rZ g\|_Z = \|g\|_{Z'}$. \om{Thus, $\cR(v, y)^2$
is computed for
a given $(v, y)\in V\times Y$ by introducing the lifted elements} 
\be
\label{ej}
e_j \coloneqq \rZ(A_jv- f_j),\quad j=0,\ldots,d,
\ee
so that, by linearity
\be
\label{Q}
\cR(v,y)^2 = \Big\| e_0 + \sum_{j=1}^d y_j e_j\Big\|^2_Z .
\ee
}

\om{Note that} the above derivation is still idealized \dwo{as the $d+1$ variational problems \eref{ej} are posed 
in the infinite dimensional space $Z$. As} already stressed in Remark \ref{spacediscretization}, all computations take place
in  \dwo{reference finite element spaces $V_h\subset V$ and $Z_h\subset Z$}. We thus approximate 
the $e_j$ by $e_{j,h}\in Z_h$, for $v\in V_h$,
using the Galerkin approximation of \eref{Riesz}.
 This gives rise to a computable least-squares functional 
 \be \label{Rvy0h}
\cR_h(v,y)^2 = \Big\| e_{0,h} + \sum_{j=1}^d y_j e_{j,h}\Big\|^2_Z ,\quad y\in Y.
\ee
The practical distance surrogate is then defined through the
corresponding constrained least-squares problem 
\be \label{Svy0h}
\cS_h(v,\cM_h):= \min_{y\in Y}\cR_h(v,y),
\ee
which can be solved by
standard optimization methods. As indicated earlier,  the recovery schemes can be interpreted as taking place in a fixed discrete setting, with $\cM$ replaced by $\cM_h$,
comprised of approximate solutions in a large finite element space $V_h\subset V$, and measuring accuracy only in this
finite dimensional setting. 
One should note though that the approach allows one to disentangle
discretization errors from recovery estimates, even with regard to the underlying continuous PDE model.
In fact,  given any target tolerance $\e_h$, using a posteriori error control in $Z$, the spaces $V_h,Z_h$ can be chosen large enough
to guarantee that
\be
\label{eh}
\big|\cR(v,y) - \cR_h(v,y) \big| \le \e_h\|v\|,\quad v\in V_h.
\ee
Accordingly,   one  has $\big|\cS_h(v,\cM_h)- \cS(v,\cM)\big|\le \e_h\|v\|$, so that recovery estimates remain
meaningful with respect to the continuous setting as long as $\e_h$ remains sufficiently dominated
by the the threshholds $\e_k, \sigma, \e_{\rm noise}$ appearing in the above results. For notational simplicity
we therefore continue working in the continuous setting.

\section{Joint parameter and state estimation}\label{sec:parest}

\subsection{An estimate for $y$}

Searching for a parameter $y\in Y$, that explains an observation $w= P_W u(y)$, is a nonlinear
inverse problem. As shown next, a quantifiable estimate for $y$ can be obtained from a state estimate $u^*(w)$ combined with a residual minimization. 

\om{For any state estimate $u^*(w)$ which we compute from $w$,} the most plausible parameter is the one \om{associated to the metric projection of $u^*(w)$ into $\cM$, that is,}
$$
y^* \in \argmin_{y\in Y}\|u(y)- u^*(w)\|.
$$
\om{Note that $y^*$ depends on $w$ but we omit the dependence in the notation in what follows.} Finding $y^*$ is a difficult task since it 
requires solving a non-convex optimization
problem. 

However, \ac{as we have 
already noticed, a near metric projection of $u^*$ to $\cM$ can be computed through a \om{simple} convex problem in the case of affine parameter dependence \eref{affine},
minimizing the residual $\cR(v,y)$ given by \eref{Rvy0}.
Our estimate for the parameter is therefore
\be
\label{bary}
y^* \in \argmin_{y\in Y} \cR(u^*,y), 
\ee
and it satisfies, in view of \eref{rR},
\be
\| u^* - u(y^*)\|\le r^{-1} \cR(u^*,y^*)\le \kappa \dist(u^*,\cM),\quad \kappa = R/r.
\ee
Hence, if we use, for instance, the state estimate $u^*(w)$ from \eqref{ustar}, we conclude by Theorem \ref{theorecovery} that $u(y^*)$ deviates from the true state $u(y)$ by
\begin{align}
\|u(y) - u(y^*) \|
&\le \| u(y)- u^*(w) \|+\|u^*(w)- u(y^*)\| \nonumber\\
&\le (1+\kappa)\|u(y)- u^*(w)\| \nonumber\\
&\le   (1+\kappa)  \delta_{\kappa\sigma}, \label{bary1}
\end{align}
where $\delta_{\kappa\sigma}$ is the benchmark quantity defined in \iref{benchapp}.  If in addition also $P_W: \cM \to W$ is injective so that $\delta_0=0$ and if $W$ and $\cM$ are favorably oriented, 
as detailed in the assumptions of Theorem \ref{thm:3.3}, one even obtains}
\be
\label{better}
\|u(y) - u(y^*) \|\le (2\mu(\cM,W)+1)\kappa\sigma.
\ee
To derive from such bounds estimates for the deviation of $y^*$ from $y$, more  information on the underlying PDE  model
is needed. For instance, for the second order parametric family of elliptic PDEs \iref{ellip} and strictly positive right hand side $f$,
it is shown in \cite{BCDPW} that the parameter-to-solution map is injective. If in addition the parameter dependent diffusion coefficient
$a(y)$ belongs to $H^1(\Omega)$,
one has a quantitative {\em inverse stability} estimate of the form
\be
\label{invstab}
\|a(y)-a(y')\|_{L^2(\Omega)}\le C\|u(y)-u(y')\|^{1/6}.
\ee
Combining this, for instance, with \iref{bary1}, yields
\be
\label{parest}
\|a(y)-a(y^*)\|_{L^2(\Omega)}\le C(1+\delta)^{1/6}\delta_{\kappa\sigma}^{1/6}.
\ee
\ac{Under the favorable assumptions of Theorem 
\ref{thm:3.3}, one obtains a bound of the form
\be
\|a(y)-a(y^*)\|_{L^2(\Omega)}\lesssim \sigma^{1/6}.
\ee
Finally, in relevant situations (Karhunen-Loeve expansions)  the functions $\psi_j$ in the expansion of $a(y)$ form an $L^2$-orthogonal system.
The above estimates translate then 
into estimates for a weighted $\ell_2$-norm,
\be
\(\sum_{j\geq 1} c_j (y_j-y_j^*)^2\)^{1/2}\lesssim \sigma^{1/6}.
\ee
where $c_j=\|\psi_j\|_{L^2}^2$.}

\subsection{Alternating residual minimization}\label{ssec:alternate}

{\anew 
The state estimate $u^*(w)$ is defined by selecting among the potential estimates $u^*_k(w)$
the one that sits closest to the solution manifold, in the sense of the surrogate distance 
$\cS(v,\cM)$. Finding the element in $V_w=w+ W^\perp$ that is closest to $\cM$ 
would provide a possibly improved state estimate, and
as pointed out in the previous section, also an improved parameter estimate. As explained earlier, it would help in addition with improved set estimators for $\cM_w$.

Adhering to the definition of the residual $\cR(v,y)$
from \iref{resid},  we are thus led to consider the double minimization problem
\be
\label{eq:min-y-gen}
\min_{(v,y) \in (w+W^\perp) \times Y} \cR(v,y)=\min_{v\in w+W^\perp} \cS(v,\cM).
\ee

We first show that a global minimizing pair $(u^*,y^*)$ 
meets the optimal benchmarks introduced in \S 2.
In the unbiased and non-noisy case, the value of the global minimum
is $0$, attained by the exact parameter $y$ and state $u(y)$. Any global minimizing pair 
$(u^*,y^*)$ will thus satisfy $P_W u^*=w$ and $u^*=u(y^*)\in \cM$. In other word, 
the state estimate $u^*$ belongs to $\cM_w$, and therefore meets the optimal
benchmark
\be
\|u-u^*\|\leq \delta_0.
\ee 
In the case of model bias and noise of amplitude, the state
$u$ satisfies
\be
{\rm dist}(u,\cM)\leq  \e_{model} \quad {\rm and} \quad \|w-P_W u\|\leq \e_{noise}.
\ee
It follows that there exists a parameter $y$ such that $\|u-u(y)\|\leq \e_{model}$
and a state $\t u\in w+W^\perp$ such that $\|u-\t u\|\leq \e_{noise}$. For this state 
and parameter, one thus have
\be
\cR(\t u, y) \leq R \|u(y)-\t u\| \leq R(\e_{model}+\e_{noise}).
\ee
Any global minimizing pair $(u^*,y^*)$ will thus satisfy 
\be
\|u^*-u(y^*)\|\leq \frac 1 r\cR(u^*, y^*) \leq \kappa(\e_{model}+\e_{noise}), 
\quad \kappa:=\frac R r.
\ee
Therefore $u^*$ belongs to the set $\cM_{\e,w}$, as defined by \iref{msigmaw},
with $\e:=\kappa(\e_{model}+\e_{noise})$ and so does $\t u$ since
$\|\t u-u(y)\|\leq \e_{model}+\e_{noise}\leq \e$. In turn, the state estimate $u^*$
meets the perturbed benchmark
\be
\|u^*-u\| \leq \e_{noise} + \|u^* -\t u\| \leq \e_{noise}+\delta_\e \leq 2\delta_{\e}.
\ee

From a numerical perspective, the search for a global minimizing pair is 
a difficult task due to the fact that $(v,y)\mapsto \cR(v,y)$ is generally not 
a convex function. However, it should be noted that in the case of 
affine parameter dependence \eref{affine},
the residual $\cR(v,y)$ given by \eref{Rvy0} is a convex function 
in each of the two variables $v,y$ separately, keeping the other one fixed. 
More precisely $(v,y)\mapsto \cR(v,y)^2$ is a quadratic convex function in each variable.
This suggests the following alternating minimization procedure.
Starting with an initial guess $u^0 \in w + W^\perp$, we iteratively compute for $j=0,1,2,..$,
\begin{align}
y^{j+1} &\in \argmin_{y\in Y} \cR(u^j,  y) \label{eq:min-y}\\
u^{j+1} &\in \argmin_{v\in V_w}\cR(v, y^{j+1}) \label{eq:min-v}.
\end{align}
Each problem has a simply computable 
minimizer, as discussed in the next section, and the residuals are non-increasing
\be
\cR(u^j,y^j)\ge \cR(u^j,y^{j+1})\ge \cR(u^{j+1},y^{j+1}) \ge \cdots
\ee
Of course, one cannot guarantee in general that $(u^j,y^j)$ 
converges to a global minimizer, and the procedure may stagnate at a local minimum. 

The above improvement property still tells us that 
if we initialize the algorithm by taking $u^0=u^*=u^*(w)$ the state estimate
from \eqref{ustar} and $y^0\in {\rm argmin}_{y\in Y}\cR(u^*,y)$, then we are
ensured at step $k$ that
\be
\cR(u^j,y^j) \leq \cR(u^*,y^*),
\ee
and therefore, by the same arguments as in the proof of Theorem \ref{thm:3.4}, 
one finds that 
\be
\|u-u^j\|\leq \delta_{\kappa\rho}+\e_{noise},
\ee
with $\kappa$ and $\rho$ as in \iref{pertest}. In other words, the new estimate
$u^j$ satisfies at least the same accuracy bound than $u^*$.
The numerical tests performed in \S 5.3 reveal that
it can be significanly \dw{more} accurate.
}

\subsection{Computational issues}

We now explain how to efficiently compute the steps in \eqref{eq:min-y} and \eqref{eq:min-v}. We continue to consider a family of linear parametric PDEs  with affine parameter dependence \iref{affine}, admitting a uniformly stable variational formulation over the pair 
trial and test spaces $V, Z$, see \iref{varform}-\iref{boundA}.

\paragraph{Minimization of \eqref{eq:min-y}:} Problem \eqref{eq:min-y-gen} requires minimizing $\cR(v,y)$ for a fixed $v\in V_w$ over $y\in Y$.
According to \iref{Q}, it amounts to solving a constrained linear least-squares problem
\be
\label{LS}
\min_{y\in Y} \Big\| e_0 + \sum_{j=1}^d y_j e_j\Big\|^2_Z ,
\ee
where the $e_j\in Z$ are the Riesz-lifts $\rZ(A_jv- f_j)$, $j=0,\ldots,d$, defined in \iref{ej}. As indicated earlier, the 
numerical solution of \iref{LS} (for $e_j=e_{j,h}\in Z_h\subset Z$) is standard.

\paragraph{Minimization of \eqref{eq:min-v}:} Problem \eqref{eq:min-v-gen} is of the form
\be
\label{eq:min-v-gen}
\min_{v\in w+W^\perp} \cR(v,  y)^2 = \min_{v\in w+W^\perp} \Vert A(y) v - f(y) \Vert^2_{Z'} 
\ee
for a fixed $y\in Y$. A naive approach for solving \eqref{eq:min-v-gen} would consist in working in a closed subspace of $\widetilde W^\perp \subseteq W^\perp$ of sufficiently large dimension. We would then optimize over $v \in w + \widetilde W^\perp$. However, this would lead to a large quadratic problem of size  $\dim \widetilde W^\perp$ which would involve $\dim \widetilde W^\perp$ Riesz representer computations. We next propose an alternative strategy involving the solution of only $m+3$ variational problems. To that end, we assume in what follows that $V$ is continuously embedded in $Z'$, which is the case for all the examples of interest, mentioned earlier in the paper.

{\anew The proposed strategy is based on
 two isomorphisms from $V$ to $Z$ that preserve inner products in a  sense
 to be explained next. We make again heavy use of the Riesz-isometry
defined in \iref{Riesz} and consider the two isomorphisms
\be
\label{TC}
T=T(y)  := \rZ A(y) : V \to Z,\quad S=S(y) := A(y)^{-*}\rV^{-1} : V \to Z,
\ee
where $r_Z:Z'\to Z$ and $r_V: V'\to V$ are the previously introduced
Riesz lifts. One then observes that, by standard duality arguments,
they preserve inner products  in the sense that for $u,v\in V$
\be
\label{TC2}
\langle Tu ,S v\rangle_Z = \langle \rZ A(y)u,A(y)^{-*}\rV^{-1}v\rangle_Z  = \langle u,v\rangle_V,
\ee
where we have used selfadjointness of Riesz isometries.
In these terms the objective functional $\cR(v,y)^{\dw{2}}$ takes the form
\be
\Vert A(y) v - f(y) \Vert^2_{Z'} = \Vert Tv - \rZ f(y)\|^2_Z.
\ee
We can use \eqref{TC2} to reformulate \eqref{eq:min-v-gen} as
\be
\min_{v\in w+W^\perp} \cR(v,  y)^2 
= \min_{v\in w+W^\perp} \Vert Tv - \rZ f(y) \Vert^2_{Z} 
= \min_{z \in Tw+ S(W)^\perp} \Vert z - \rZ f(y) \Vert^2_{Z} ,
\ee
\om{where we have used that $T(W^\perp)=S(W)^\perp$ to obtain the last equality.} Note that the unique solution $z^*\in Z$ to the right hand side gives a solution $v^*\in V$ to the original problem through the relationship $T v^* = z^*$.
The minimizer $z^*$ can be obtained by an appropriate orthogonal projection onto $S(W)$. This
indeed amounts to solving a fixed number of $m+3$ variational problems
without compromising accuracy by choosing a perhaps too moderate
dimension for a subspace $\wt W^\perp$ of $W^\perp$.

More precisely, we have $z^*=Tw+\t z$ where
$\t z\in S(W)^\perp$ minimizes $\Vert \t z+Tw - \rZ f(y) \Vert^2_{Z}$, and therefore
\be
\t z=P_{S(W)^\perp}(\rZ f(y)-Tw)=
\rZ f(y)-Tw-P_{S(W)}(\rZ f(y)-Tw).
\ee
This shows that
\be
\label{z*1}
z^* =z^*(y) :=  f(y) - P_{S(W)} (\rZ f(y)-Tw).
\ee
Thus, a single iteration of the type \iref{eq:min-v-gen} requires assembling $z^*$ followed by solving the variational problem
\be
\label{PG}
\langle Tv^*,z\rangle_Z = (A(y)v^*)(z)= \langle z^*,z\rangle_Z,\quad z\in Z,
\ee
that gives $v^*$. Assembling $z^*$ involves
\begin{enumerate}
\item[(i)] evaluating $Tw$, which means solving the Riesz-lift $\langle Tw,z\rangle_Z = (A(y)w)(z)$, $z\in Z$;
\item[(ii)] computing the Riesz-lift $\rZ f(y)$ by solving $\langle \rZ f(y),z\rangle_Z = (f(y))(z)$, $z\in Z$;
\item[(iii)] computing the projection $P_{S(W)} (\rZ f(y)-Tw)$. This requires computing the transformed basis functions 
$Sw_i=A(y)^{-*}\rV^{-1}w_i$,
which are solutions to the variational problems
\be
\label{Cwi}
 (A(y)^* Sw_i)(v)= \langle w_i,v\rangle_V,\,\,v\in V, \quad i=1,\ldots,m.
\ee
\end{enumerate}
}

Of course, these variational problems are solved only approximately in appropriate large but finite dimensional
spaces $V_h\subset V, Z_h\subset Z$ along the remarks at the end of the previous section. While approximate 
Riesz-lifts involve symmetric variational formulations which are well-treated by Galerkin schemes, the 
problems involving the operator $A(y)$ or $A(y)^*$ may in general require an unsymmetric variational formulation
where $Z\neq V$ and Petrov-Galerkin schemes on the discrete level. 
For each of the examples (such as a time-space variational formulation of parabolic or 
convection diffusion equations) stable discretizations are known, see e.g. \cite{BS,CDW,DHSW,DPW,SW}.

A particularly important strategy for unsymmetric problems is to write the PDE first as an equivalent system of first order PDEs
permitting a so called ``ultraweak'' formulation where the (infinite-dimensional) trial space $V$ is actually
an $L_2$-space and the required continuous embedding $V\subset Z'$ holds. 
The mapping $\rV$ is then just the identity and so called Discontinuous Petrov Galerkin
methods offer a way of systematically finding appropriate test spaces in the discrete case with 
uniform inf-sup stability, \cite{CaDeGo}. In this context, the mapping $T$ from\iref{TC} plays a pivotal role 
in the identification of ``optimal test spaces'' and is referred to as ``trial-to-test-map''. 

Of course, in the case of problems that admit a symmetric variational formulation, i.e., $V=Z$, things simplify
even further. To exemplify this, consider the a parametric family of elliptic PDEs \iref{ellip}. In this case one has
(assuming homogeneous boundary conditions) $V=Z=H^1_0(\Omega)$ so that $\rZ=\rV = \Delta^{-1}$.
Due to the selfadjointness of the underlying elliptic operators $A(y)$ in this case, the problems \iref{Cwi}
are of the same form as in \iref{PG} that can be treated on the discrete level by standard Galerkin discretizations.

\section{Numerical illustration}

In this section we illustrate the construction of nonlinear reduced models, and demonstrate the mechanism of model selection using the residual surrogate methods outlined in \S \ref{ssec:residual}.

In our tests we consider the elliptic problem mentioned in \S \ref{ssec:obj} on the unit square $D=]0,1[^2$ with homogeneous Dirichelet boundary conditions, and a parameter dependence in the diffusivity field $a$. Specifically, we consider the problem
\be \label{diff_eq_y}
-{\rm div}(a(y) \nabla u)=f, 
\ee
with $f=1$ on $D$, with $u_{|\partial D}=0$. The classical variational formulation uses the same trial and test space 
$V=Z=H^1_0(D)$. We perform space discretization by the Galerkin 
method using $\mathbb{P}_1$ finite elements to produce solutions $u_h(y)$, with a triangulation on a regular grid of mesh size $h=2^{-7}$. 

\subsection{Test 1: pre-determined splittings}

\begin{figure}
\centering
\begin{tikzpicture}
\draw [thick] (0,0) rectangle (4,4);
\draw (0,0) rectangle (3,3);
\draw (3,3) rectangle (4,4);

\node at (1.5,1.5) {$D_{1,1}$};
\node at (1.5,3.5) {$D_{1,2}$};
\node at (3.5,1.5) {$D_{1,3}$};
\node at (3.5,3.5) {$D_{1,4}$};

\draw[->] (0,0) -- (4.6,0) node[anchor=west] {$x_1$};
\draw[->] (0,0) -- (0,4.6) node[anchor=south] {$x_2$};
\node [below left] at (0,0) {$0$};

\draw [thick] (6,0) rectangle (10,4);
\draw (6,0) rectangle (7,1);
\draw (7,1) rectangle (10,4);

\draw[->] (6,0) -- (10.6,0) node[anchor=west] {$x_1$};
\draw[->] (6,0) -- (6,4.6) node[anchor=south] {$x_2$};
\node [below left] at (6,0) {$0$};

\node at (6.5,0.5) {$D_{2,4}$};
\node at (6.5,2.5) {$D_{2,3}$};
\node at (8.5,0.5) {$D_{2,2}$};
\node at (8.5,2.5) {$D_{2,1}$};
\end{tikzpicture}
\caption{\footnotesize Left, the partition of the unit square used in $a_1$, and right, the partition used in $a_2$.}
\label{fig:partitions}
\end{figure}
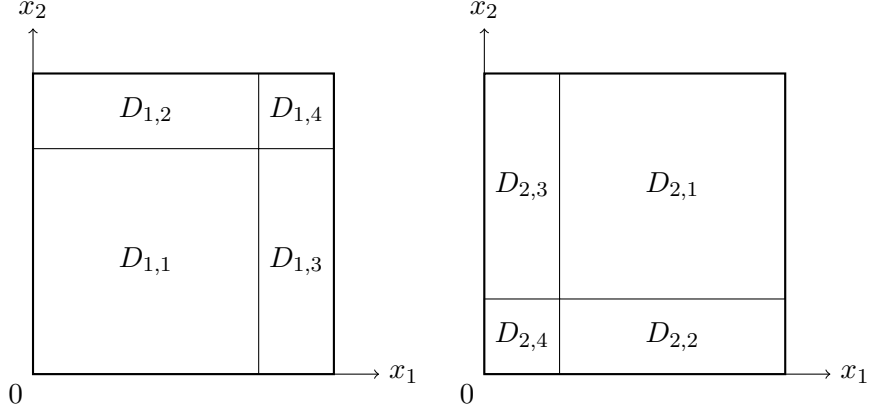

In this first test, we examine the reconstruction performance with localized reduced bases on a manifold having a predetermined splitting.
Specifically, we consider two partitions of the unit square, $\{ D_{1, \ell} \}_{\ell=1}^4$ and $\{ D_{2,\ell} \}_{\ell=1}^4$, with
\begin{alignat*}{7}
&D_{1,1} \coloneqq \[0,\frac 3 4\[ \times \[0, \frac 3 4\[& \quad &D_{1,2} \coloneqq \[0,\frac 3 4\[ \times \[\frac 3 4, 1\]& \quad &D_{1,3} \coloneqq \[\frac 3 4,1\] \times \[0, \frac 3 4\[& \quad &D_{1,4} \coloneqq \[\frac 3 4, 1\] \times \[ \frac 3 4, 1\]& \\
&D_{2,1} \coloneqq \[\frac 1 4, 1\] \times \[\frac 1 4, 1\]& \quad &D_{2,2} \coloneqq \[\frac 1 4,1\] \times \[0, \frac 1 4\[& \quad &D_{2,3} \coloneqq \[0,\frac 1 4\[ \times \[\frac 1 4,1\]& \quad &D_{2,4} \coloneqq \[0,\frac 1 4\[ \times \[0, \frac 1 4\[&
\end{alignat*}
The partitions are symmetric along the axis $x+y = 1$ as illustrated in Figure \ref{fig:partitions}. This will play a role in the understanding of the results below. We next define two parametric diffusivity fields
\be \label{def:a_1_a_2}
a_1(y) \coloneqq \overline a + c \sum_{\ell = 1} ^4 \Chi_{D_{1,\ell}} y_\ell \quad \text{and} \quad 
a_2(y) \coloneqq \overline a + c \sum_{\ell = 1} ^4 \Chi_{D_{2,\ell}} y_\ell \, ,
\ee
where the vector of parameters $y=(y_1,\dots,y_4)$ ranges in $Y = [-1, 1]^4$ and $\Chi_{D_{1,\ell}}$ is the indicator function of $D_{1,\ell}$ (similarly for $\Chi_{D_{2,\ell}}$). The fields $a_1(y)$ and $a_2(y)$ are mirror images of each other along $x+y = 1$. In the numerical tests that follow, we take $\overline a = 1$ and $c = 0.9$.

We denote by $u_1(y)$ the solution to the elliptic problem \eqref{diff_eq_y} with diffusivity field $a_1(y)$, and then label by 
$\cM_{1}:=\{u_1(y)\;:\; y\in Y\}$ the resulting solution manifold. Strictly speaking, we should write $\cM_{h,1}$ as our solutions are finite dimensional approximations, however we suppress the $h$ as there should be little ambiguity going forward. Similarly, $\cM_2$ is the set of all solutions $u_2(y)$ of \eqref{diff_eq_y} over $Y$ where the diffusivity field is given by $a_2$. We take their union $\cM = \cM_1 \cup \cM_2$ to be our global solution manifold that has the obvious pre-determined splitting available to us. 

For our computations, we generate training and test sets. For the training, we draw $N_{\mathrm{tr}}=5000$ independent samples $\wt Y_{\mathrm{tr}}=( y^{\mathrm{tr}}_j)_{j=1}^{N_\mathrm{tr}}$ that are uniformly distributed over $Y$. The collection of solutions $\wt \cM_1 \coloneqq \{ u_1(y^{\mathrm{tr}}_j) \}_{j=1}^{N_\mathrm{tr}}$ and $\wt \cM_2 \coloneqq \{ u_2(y^{\mathrm{tr}}_j) \}_{j=1}^{N_\mathrm{tr}}$, are used as {\em training sets} for $\cM_1$ and $\cM_2$. The training set for the full manifold $\cM$ is $\wt\cM = \wt\cM_1 \cup \wt\cM_2$. Since we use the same parameter points $y^{\mathrm{tr}}_j$ for both sets, any solution in $\wt \cM_1$ has a corresponding solution in $\wt \cM_2$ that is its symmetric image along the axis $x+y=1$. To test our reconstruction methods, we generate $N_{\mathrm{te}}=2000$ independent points in $Y$ that are distinct from the training set. The corresponding test solution sets are $\wt \cT_1$ and $\wt \cT_2$. All computations are done by solving \eqref{diff_eq_y} in the finite element space.

Given an unknown $u\in\cM$, we want to recover it from its observation $w=P_Wu$. For the measurement space $W$, we take a collection of $m=8$ measurement functionals $\ell_i(u) = \<\omega_i, u\> = |B_i|^{-1} \int u \, \Chi_{B_i}$ that are local averages in a small area $B_i$ which are boxes of width $2h = 2^{-6}$, each placed randomly in the unit square. The measurement space is then $W = \Span\{\omega_1,\ldots, \omega_m\}$.

\begin{figure}
\centering   
  \includegraphics[width=\linewidth]{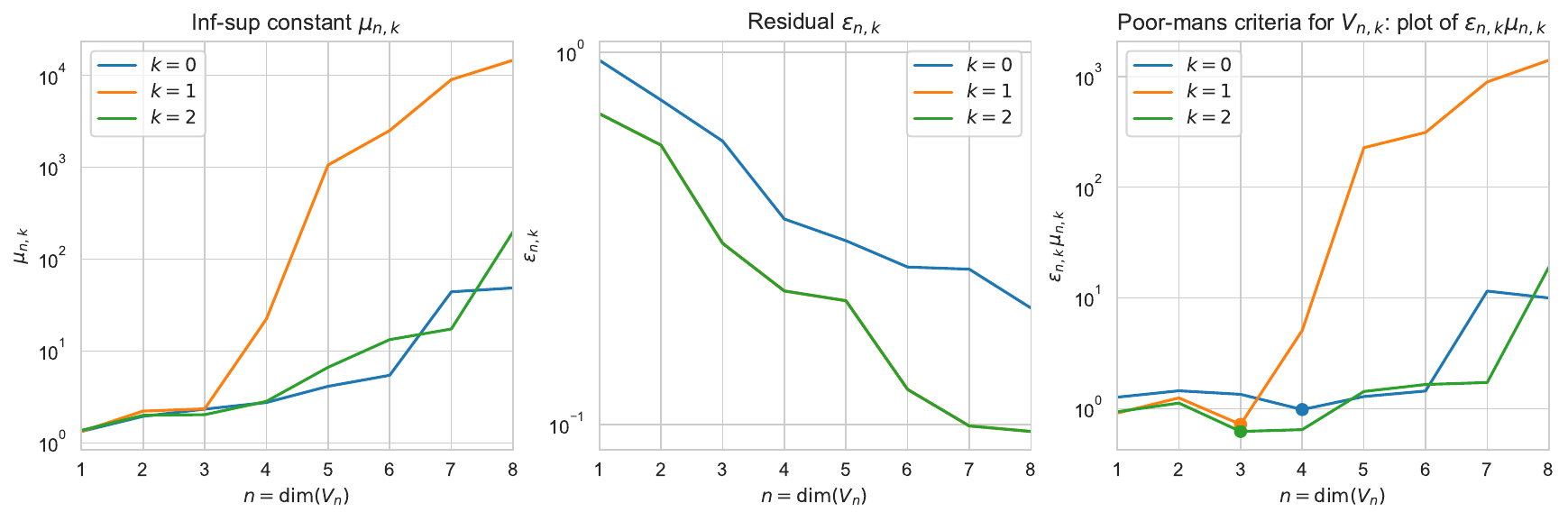}
  \vspace{-20pt}
  \caption{\footnotesize Inf-sup constants $\mu_{n,k}$ and residual errors $\e_{n,k}$, leading to the dimension $n$ choice for $V_k$.}
  \label{fig:t1_mu_eps}
\end{figure}

Since we are only given $w$, we do not know if the function to reconstruct is in $\cM_1$ or $\cM_2$ and we consider two possibilities for reconstruction:
\begin{itemize}
\item \emph{Affine method:} We use affine reduced models $V_{n,0} = \bar u_0 + \bar V_{n,0}$ generated for the full manifold $\cM = \cM_1 \cup \cM_2$. In our example, we take $\bar u_0 = u(y=0)$ and $\bar V_{n,0}$ is computed by
the greedy selection algorithm over $\wt\cM-\bar u_0$.
Of course the spaces $V_{n,0} $ with $n$ sufficiently large have high potential for approximation of the full manifold $\cM$, and obviously also for the subsets $\cM_1$ and $\cM_2$ (see Figure \ref{fig:t1_mu_eps}). Yet, we can expect some bad artefacts in the reconstruction with this space since the true solution will be approximated by snapshots, some of which coming from the wrong part of the manifold and thus associated to the wrong partition of $D$. In addition, we can only work with $n\leq m=8$ and this may not be sufficient regarding the approximation power. Our estimator $u^*_0(w)$ uses the space $V_{n^*_0,0}$, where the dimension $n^*_0$ is the one that reaches $\tau_0 = \min_{1\leq n \leq m} \mu_{n,0}\e_{n,0}$ as defined in \eqref{eq:acc_and_stab} and \eqref{eq:pm_criteria}. Figure \ref{fig:t1_mu_eps} shows the product $\mu_{n,0}\e_{n,0}$ for $n=1,\dots,m$ and we see that $n^*_0=3$.

\item \emph{Nonlinear method:} We generate affine reduced bases spaces $V_{n,1}=\bar u_1 + \bar V_{n,1}$ and $V_{n,2}=\bar u_2 + \bar V_{n,2}$, each one specific for $\cM_1$ and $\cM_2$. Similarly as for the affine method, we take as offsets $\bar u_i=u_i(y=0)=\bar u_0$, for $i=1,2$, and we run two separate greedy algorithms over $\wt\cM_1-\bar u_1$ and $\wt\cM_2-\bar u_2$ to build $\bar V_{n,1}$ and $\bar V_{n,2}$. We select the dimensions $n^*_k$ that reach $\tau_{k}=\min_{n=1,\dots, m}\mu_{n,k} \e_{n,k}$ for $k=1,2$. From Figure \ref{fig:t1_mu_eps}, we deduce\footnote{Due to the spatial symmetry along the axis $x+y=1$ for the functions in $\wt\cM_1$ and $\wt\cM_2$, the greedy algorithm selects exactly the same candidates to build $V_{n,2}$ as for $V_{n,1}$, except that each element is mirrored in the axis. One may thus wonder why $n^*_1\neq n^*_2$. The fact that different values are chosen for each manifold reflects the fact that the measurement space $W$ introduces a symmetry break and the reconstruction scheme is no longer spatially symmetric contrary to the $V_{n,k}$.} that $n^*_1=4$ and $n^*_2=3$. This yields two estimators  $u^*_1(w)$ and $u^*_2(w)$. We can expect better results than the affine approach if we can detect well in which part of the manifold the target function is located. The main question is thus whether our model selection strategy outlined in Section \ref{sec:model-selection} is able to detect well from the observed data $w$ if the true 
$u$ lies in $\cM_1$ or $\cM_2$. For this, we compute the surrogate manifold distances
\begin{align}
\cS(u^*_k(w), \cM_k) \coloneqq \min_{y\in Y} \cR_k( u^*_k(w), y),\quad k=1,2,
\label{eq:Sk}
\end{align}
where
$$
\cR_k( u^*_k(w), y) \coloneqq \| \mathrm{div}(a_k(y) \nabla u^*_k(w)) \dw{+} f \|_{V'}
$$
is the residual of $u^*_k(w)$ related to the PDE with diffusion field $a_k(y)$. To solve problem \eqref{eq:Sk}, we follow the steps given in Section \ref{ssec:residual}. The final estimator is $u^*=u^*_{k^*}$, where
$$
k^* = \argmin_{k=1,2} \cS(u^*_k(w), \cM_k).
$$
\end{itemize}

\begin{figure}
\centering   
  \includegraphics[width=\linewidth]{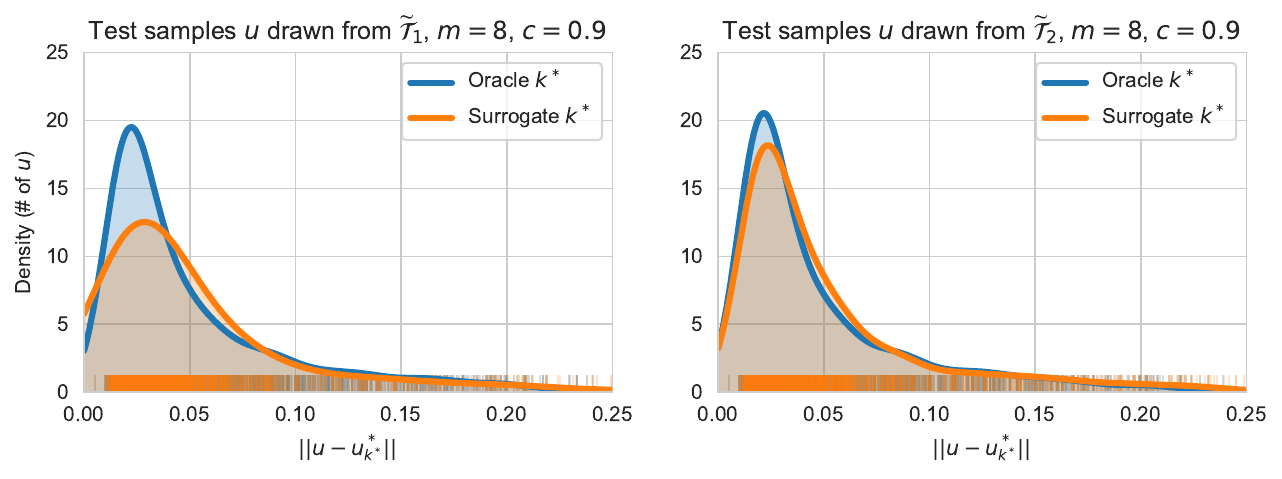}
  \vspace{-20pt}
  \caption{\footnotesize Kernel density estimate (KDE) plot of the $u_1^*$ and $u_2^*$.}
  \vspace{-4pt}
  \label{fig:t1_err_densities}
\end{figure}

Table \ref{tab:t1_counts} quantifies the quality of the model selection approach. It displays how many times our model selection strategy yields the correct result $k^*=1$ or incorrect result $k^*=2$ for the functions from the test set $\wt\cT_1\subset \cM_1$ (and vice-versa for $\wt\cT_2$). Recalling that these tests sets have $N_{te}=2000$ snapshots, we conclude that the residual gives us the correct manifold portion roughly 75\% of the time. We can compare this performance with the one given by the oracle estimator (see Table \ref{tab:t1_counts})
$$
k^*_{\text{oracle}}=\argmin_{k=1,2} \| u - u^*_k(w) \| .
$$
In this case, we see that the oracle selection is very efficient since it gives us the correct manifold portion roughly 99\% of the time. Figure \ref{fig:t1_err_densities} completes the information given in Table \ref{tab:t1_counts}  by showing the distribution of the values of the residuals and oracle errors. The distributions give visual confirmation that both the model and oracle selection tend to pick the correct model by giving residual/error values which are lower in the right manifold portion. Last but not least, Figure \ref{fig:t1_mu_eps} gives information on the value of inf-sup constants and residual errors leading to the choice of the dimension $n^*$ for the reduced models. Table \ref{tab:errors-ex1} summarizes the reconstruction errors.


\begin{table}
\centering\footnotesize
\begin{minipage}{.49\linewidth}
\begin{tabular}{| p{2.9cm} | p{1.8cm} | p{1.8cm} |}
  \multicolumn{1}{r}{} & \multicolumn{2}{ c }{Test function from:} \\ \hline
\textbf{Surrogate selection} & $\wt\cT_1\subset \cM_1$ & $\wt\cT_2\subset \cM_2$ \\
 \hline
 $k^*=1$ & \cellcolor{green!25}1625 & \cellcolor{red!25}386\\
 \hline
  $k^*=2$ & \cellcolor{red!25}375 & \cellcolor{green!25}1614\\
 \hline
Success rate & 81.2 \% & 80.7 \% \\
 \hline
\end{tabular}
\end{minipage}
\begin{minipage}{.49\linewidth}
\begin{tabular}{| p{2.9cm} | p{1.8cm} | p{1.8cm} |}
  \multicolumn{1}{r}{} & \multicolumn{2}{ c }{Test function from:} \\ \hline
\textbf{Oracle selection} & $\wt\cT_1\subset \cM_1$ &  $\wt\cT_2\subset \cM_2$ \\
 \hline
 $k^*=1$ & \cellcolor{green!25}1962 & \cellcolor{red!25} 9\\
 \hline
  $k^*=2$ & \cellcolor{red!25}38 & \cellcolor{green!25}1991\\
 \hline
Success rate & 98.1 \% & 99.5 \% \\
 \hline
\end{tabular}
\end{minipage}
\caption{Performance of model selection and oracle selection. \label{tab:t1_counts}}
\end{table}

\begin{table}
\centering \footnotesize
\begin{minipage}{.49\linewidth}
\begin{tabular}{| p{3.3cm} | p{1.6cm} | p{1.6cm} |}
 \multicolumn{1}{r}{} & \multicolumn{2}{ c }{Test function from:} \\ \hline
\textbf{Average error } & $\wt\cT_1\subset \cM_1$ & $\wt\cT_2\subset \cM_2$ \\
 \hline
 Affine method $u^*_0$ & 6.047e-02   & 6.661e-02 \\
 \hline
  Nonlinear with oracle model selection & 5.057e-02 & 4.855e-02 \\
 \hline
 Nonlinear with surrogate model selection & 5.522e-02 & 5.201e-02 \\
 \hline
\end{tabular}
\end{minipage}
\begin{minipage}{.49\linewidth}
\begin{tabular}{| p{3.3cm} | p{1.6cm} | p{1.6cm} |}
 \multicolumn{1}{r}{} & \multicolumn{2}{ c }{Test function from:} \\ \hline
\textbf{Worst case error } & $\wt\cT_1\subset \cM_1$ & $\wt\cT_2\subset \cM_2$ \\
 \hline
 Affine method $u^*_0$ & 4.203e-01 & 4.319e-01 \\
 \hline
  Nonlinear with oracle model selection & 2.786e-01 & 2.641e-01 \\
 \hline
 Nonlinear with surrogate model selection & 4.798e-01  & 2.660e-01 \\
 \hline
\end{tabular}
\end{minipage}
\caption{Reconstruction errors with the different methods. \label{tab:errors-ex1}}
\end{table}


\newpage
\subsection{Test 2: constructing $\sigma$-admissible families} \label{sec:test_2}

\om{In this example we examine the behavior of the splitting scheme to construct $\sigma$-admissible families outlined in \S \ref{ssec:admissible_families}.}

The manifold $\cM$ is given by the solutions to equation \eqref{diff_eq_y} associated to the diffusivity field
\be \label{def_a_test}
a(y) = \overline a + \sum_{\ell = 1}^{d} c_\ell \Chi_{D_\ell} y_\ell,\quad y \in Y,
\ee
where $\Chi_{D_\ell}$ is the indicator function on the set $D_\ell$, and parameters ranging uniformly in $Y=[-1,1]^d$.
We study the impact of the intrinsic dimensionality of the manifold by considering two cases for the partition of the unit square $D$, a $2\times 2$ uniform grid partition resulting in $d=4$ parameters, and a $4\times 4$ grid partition of $D$ resulting in $d=16$ parameters. 
We also study the impact of coercivity and anisotropy on our reconstruction algorithm by examining the different manifolds generated by taking $c_\ell = c_1 \ell^{-r}$
with $c_1=0.9$ or $0.99$ and $r=1$ or $2$. The value $c_1=0.99$ corresponds to
a severe \dw{degeneration} of coercivity, and the rate $r=2$ corresponds to a more pronounced
anisotropy.

%

\begin{figure} 
\centering   
  \includegraphics[width=0.6\linewidth]{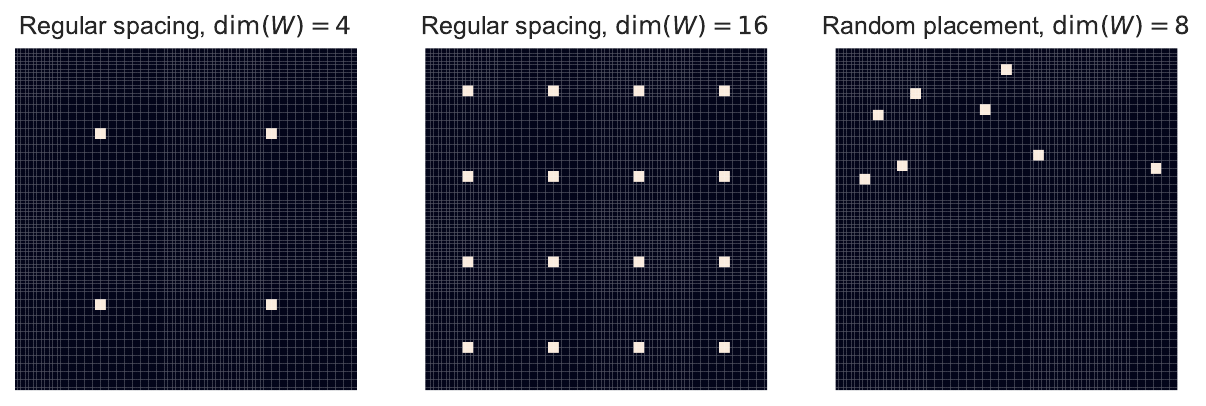}
  \vspace{-4pt}
  \caption{Measurement locations for Test 2 and Test 3}
\vspace{-8pt}
\label{fig:meas_placement}
 \end{figure}
  
We use two different measurement spaces, one with $m=\dim(W)=4$ evenly spaced local averages and the other with $m=16$ evenly spaced local averages. The measurement locations are shown diagrammatically in Figure \ref{fig:meas_placement}. The local averages are exactly as in the last test, taken in squares of side-length $2^{-6}$. Note that the two values $m=4$ and $m=16$ which we consider for the dimension of the measurement space are the same as the parameter dimensions $d=4$ and $d=16$ of the manifolds. This allows \jn{us} to study different regimes:
\begin{itemize}
\item When $m<d$, we have a highly ill-posed problem since the intrinsic dimension of the manifold is larger than the dimension of the measurement space. In particular, we expect that the fundamental barrier $\delta_0(\cM)$ is strictly positive.
Thus we cannot expect very accurate reconstructions even with the splitting strategy.
\item When $m \geq d$, the situation is more favorable and we can expect that the reconstruction involving manifold splitting brings significant accuracy gains. 
\end{itemize}

As in the previous case, the training set $\wt\cM$ is generated by a subset  $\wt Y_{\mathrm{tr}}=\{ y^{\mathrm{tr}}_j \}_{j=1,\ldots,N_\mathrm{tr}} $ of $N_{\mathrm{tr}}=5000$ samples taken uniformly on $Y$. We build the $\sigma$-admissible families outlined in \S \ref{ssec:admissible_families} using a dyadic splitting and the splitting rule is given by \eqref{eq:pm_criteria}. For example, our first split of $Y$ results in two rectangular cells $Y_1$ and $Y_2$, and the corresponding collections of parameter points $\wt Y_1 \subset Y_1$ and $\wt Y_2 \subset Y_2$, as well as split collections of solutions $\wt \cM_1$ and $\wt \cM_2$. On each $\wt \cM_k$ we apply the greedy selection procedure, resulting in $V_k$, with computable values $\mu_k$
and $\e_k$. The coordinate direction in which we split $Y$ is precisely the direction that gives us the smallest resulting $\sigma = \max_{k=1,2} \mu_k \e_k$, so we need to compute greedy reduced bases for each possible splitting direction before deciding which results in the lowest $\sigma$. Subsequent splittings are performed in the same manner, but at each step we first chose cell  $k_{\mathrm{split}} = \argmax_{k=1,\ldots,K} \mu_k \e_k$ to be split.

After $K-1$ splits, the parameter domain is divided into $Y = \bigcup_{k=1}^K Y_k$ disjoint subsets $Y_k$ and we have \dw{computed} a family of $K$ affine reduced spaces $(V_k)_{k=1,\ldots,K}$. For a given $w\in W$, we have $K$ possible reconstructions $u_1^*(w),\dots u_K^*(w)$ and we select a value $k^*$ with the surrogate based model selection outlined in \S \ref{ssec:residual}. The test is done on a test set of $N_{\text{te}}=1000$ snapshots which are different from the ones used for the training set $\wt\cM$.

\begin{figure}[H]
\centering
  \includegraphics[width=0.95\linewidth]{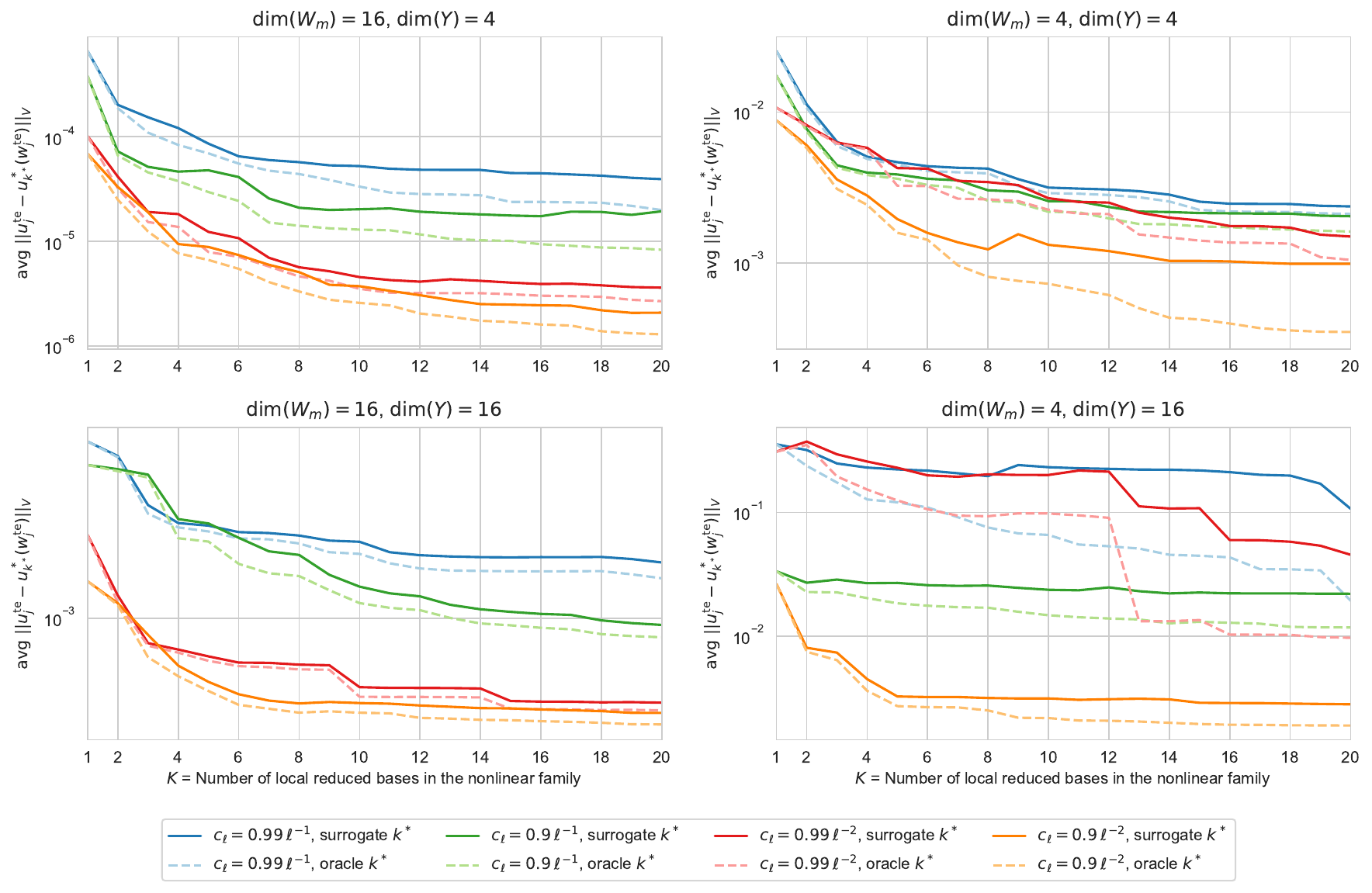}
  \vspace{-8pt}
  \caption{\footnotesize Average of errors $\| u^{\mathrm{te}}_j - u^*_{k^*}(w^{\mathrm{te}}_j ) \| $ for different choices of $k^*$ .}
  \vspace{-8pt}
  \label{fig:t2_split_errs}
\end{figure}

In Figure \ref{fig:t2_split_errs} we plot the reconstruction error, averaged
over the test set, as a function of the number of splits $K$ for all the different configurations: we consider the 2 different diffusivity fields $a(y)$ with $d=4$ and $d=16$ parameters, 
the two measurement spaces of dimension $m=4$ and $m=16$,
and the $4$ different ellipticity/coercivity regimes of $c_\ell$ in $a(y)$. We also plot the error when taking for $k^*$ the {\em oracle value} that corresponds to the value 
of $k$ that contains the parameter $y$ which gave rise to 
the snapshot and measurement. 

Our main findings can be summarized as follows:
\begin{enumerate}
\item The error decreases with the number of splits. As anticipated, 
the splitting strategy is more effective in the overdetermined regime $m\geq d$.
\item \dw{Degrading} coercivity has a negative effect on the estimation
error, while anisotropy has a positive effect.
\item Choosing $k^*$ by the surrogate based model selection yields error curves
that are above yet close to those obtained with the oracle choice. The largest
discrepancy is observed when coercivity degrades.
\end{enumerate}

\begin{figure}[H]
\centering
  \includegraphics[width=0.95\linewidth]{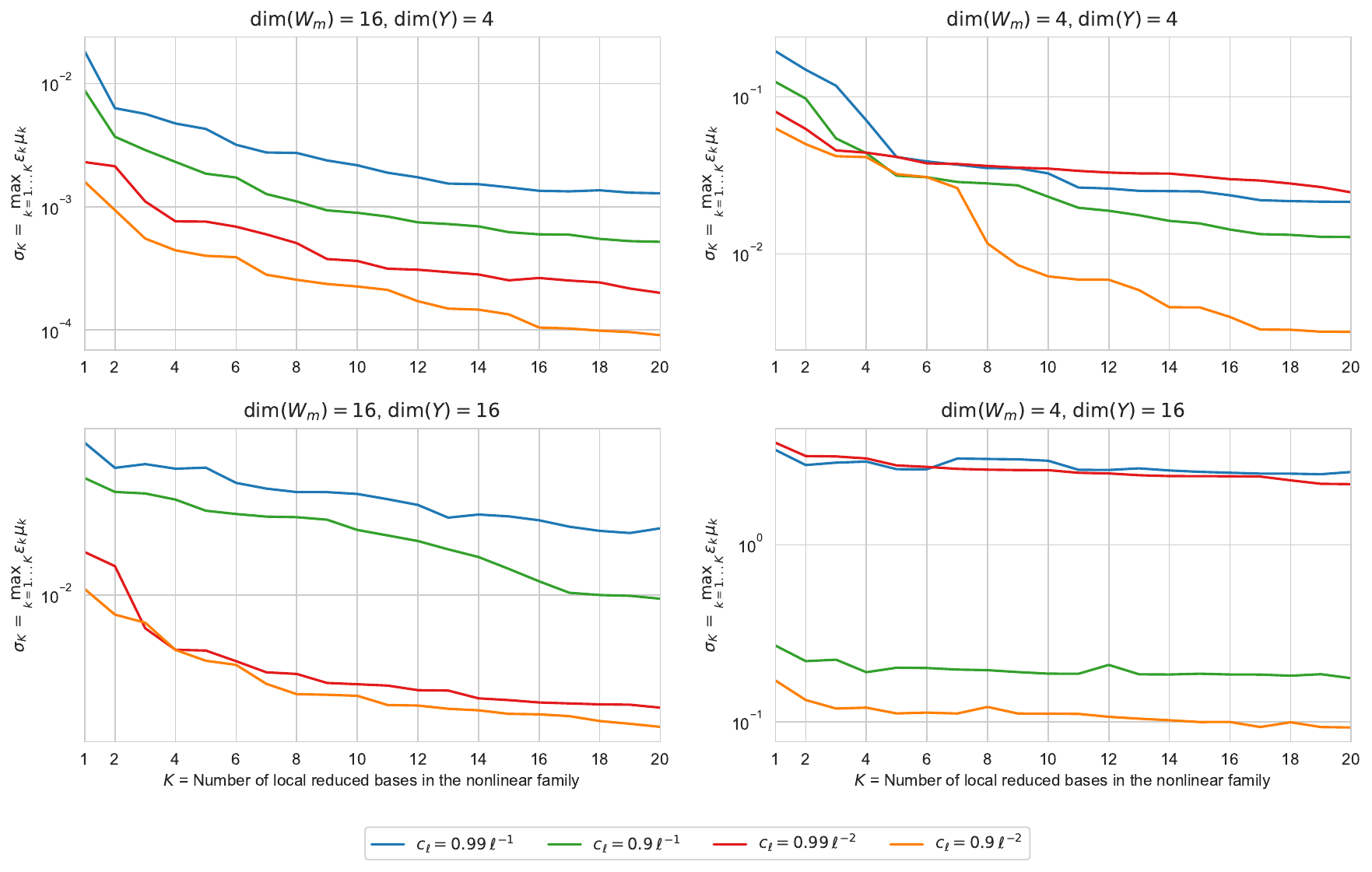}
  \vspace{-8pt}
  \caption{\footnotesize Error bounds of local linear families, given by $\sigma_K = \max_{k=1\ldots K} \mu_k \e_k$.}
  \vspace{-8pt}
  \label{fig:t2_split_sigma}
\end{figure}

\begin{figure}[H]
\centering   
  \includegraphics[width=0.95\linewidth]{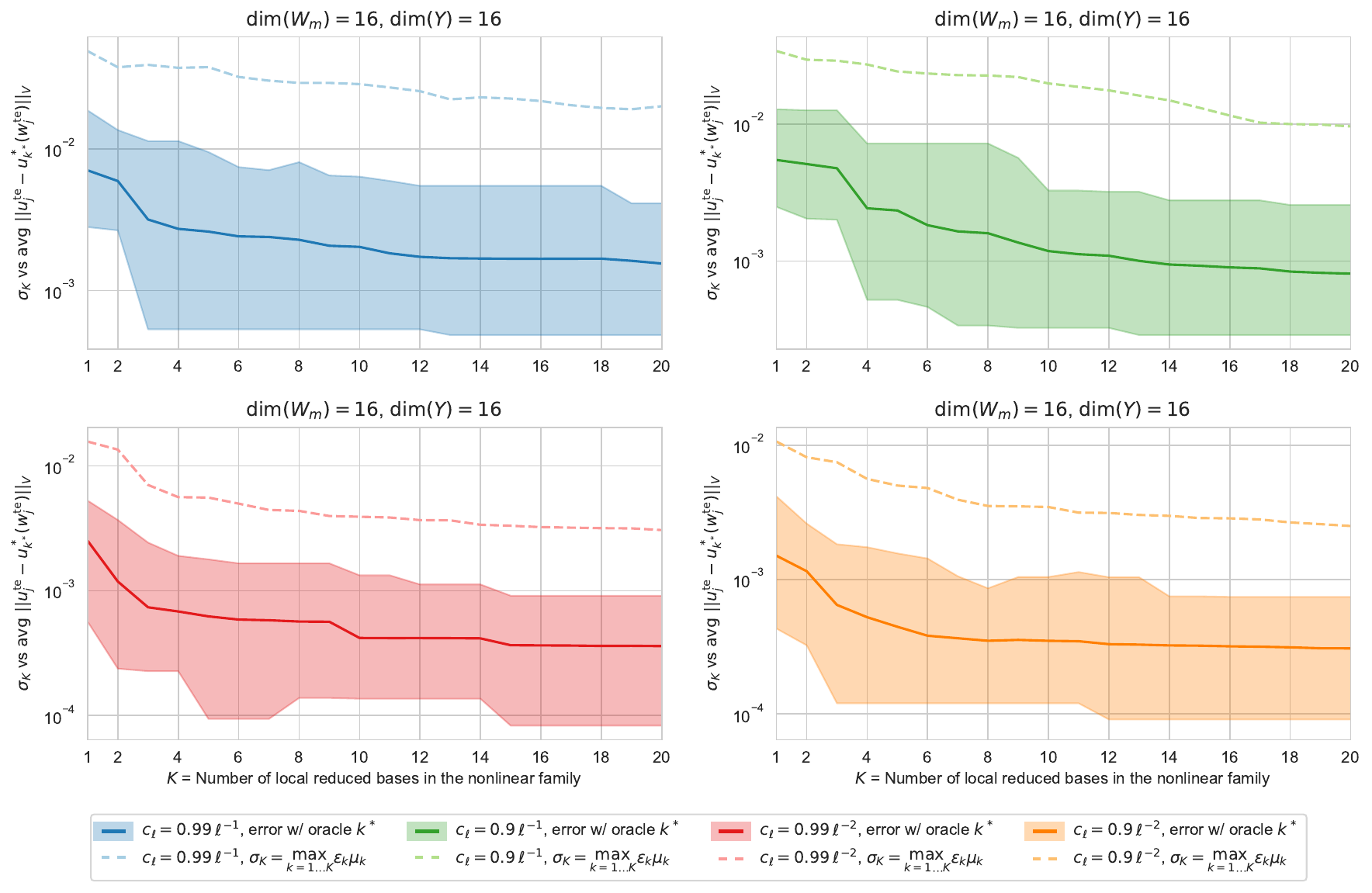}
  \vspace{-8pt}
  \caption{\footnotesize Comparison between $\sigma_K$ (dashed curve),
  the averaged oracle error (full curve) and the 
  the range from maximum to minimum oracle error (shaded region).}
  \vspace{-8pt}
  \label{fig:t2_split_sigma_vs_errs}
\end{figure}

Figure \ref{fig:t2_split_sigma} presents the 
error bounds $\sigma_K:=\max_{k=1,\dots,K}\mu_k\e_k$ which are known
to be upper bounds for the estimation error when choosing the oracle
value for $k^*$ at the given step $K$ of the splitting procedure.
We observe that these worst upper bounds have similar behaviour
as the averaged error curves depicted on Figure \ref{fig:t2_split_errs}.
In Figure \ref{fig:t2_split_sigma_vs_errs}, for the particular
configuration $\dim(Y)=\dim(W)=16$, we demonstrate that 
$\sigma_K$ indeed acts as an upper bound for
the worst case error of the oracle estimator.

\newpage
\subsection{Test 3: improving the state estimate by alternating residual minimization}

The goal of this test is to illustrate how the alternating residual minimization
outlined in \S \ref{ssec:alternate} allows one to improve the accuracy of the state estimate.
We use the same setting as Test 2, in particular, we consider 
the solution manifold $\cM$ of equation \eqref{diff_eq_y}
with the random field $a(y)$ defined as in \eqref{def_a_test}. Again we consider the cases where the $D_\ell$ from \eqref{def_a_test} are from the $2 \times 2$ and $4 \times 4$ grid, resulting in $d=\dim(Y)=4$ and $16$ respectively. Our test uses all three measurement regimes presented in Figure \ref{fig:meas_placement}, with $m=4$ and $16$ evenly spaced local average functions, and $m=8$ randomly placed local averages confined to the upper half. We use the coercivity/anisotropy regime $c_{\ell}=0.9\, \ell^{-1}$.
 
In this test we compare three different candidates for $u^0$, the starting point of the alternating minimization procedure:
\begin{itemize}
\item 
$u^0 = w$, the measurement vector without any further approximation, or equivalently the reconstruction of minimal $H^1_0$ norm among all functions
that agree with the observations.
\item
$u^0 = u^*(w)$, the PBDW state estimation using the greedy basis 
over the whole manifold, thus starting the minimization from a ``lifted'' candidate that we hope is closer to the manifold $\cM$ and should thus offer better performance. 
\item
$u^0 = u^*_{k^*}(w)$, the surrogate-chosen local linear reconstruction from the same family of local linear models from \S \ref{sec:test_2} (where $k^*$ is the index of the chosen local linear model). In this last case we take $K=20$ local linear models, i.e.~where $Y$ has been split \dw{19} times. 
\end{itemize}

Furthermore, in the third case, we restrict our parameter range to be the local parameter range chosen by the surrogate, that is we alter the step outlined in \eqref{eq:min-y} to be
$$
y^{j+1} = \argmin_{y \in Y_{k^*}} \cR(u^{j},y),
$$
where $y^{j+1}$ denotes the parameter found at the $(j+1)$-th step of the procedure. The hope is that we have correctly chosen the local linear model and restricted parameter range from which the true solution comes thanks to our local model selection. The alternating minimization will thus have a better starting position and then a faster convergence rate due to the restricted parameter range.

In our test we use the same training set $\wt \cM$ as \dw{in} the previous test, with $N_{\text{tr}}=5000$ samples,
in order to generate the reduced basis spaces. We consider a set of $N_{\text{te}}=10$ snapshots, distinct from any snapshots in $\wt\cM$, and perform the alternate minimization for each of the snapshots in the test set.

The final figures display the state error trajectories $j\mapsto \|u-u^j\|$ for each
snapshots (dashed lines) as well as their geometric average (full lines), in different colors
depending on the initialization choice. Similarly we display the 
residual trajectories $j\mapsto \cR(u^j,y^j)= \|A(y^j)u^j-f\|_{V'}$ and 
parameter error trajectories $j\mapsto |y-y^j|_2$. Our main findings can be summarized as follows:
\begin{enumerate}
\item 
In all cases, there is a substantial gain in taking $u^0 = u^*_{k^*}(w)$, the surrogate-chosen local linear reconstruction, as starting point. In certain cases, the iterative procedure initiated from the two other choices $w$ or $u^*(w)$ stagnates
at an error level that is even higher than $\|u-u^*_{k^*}(w)\|$. 
\item
The state error, residual and parameter error decrease to zero in the 
overdetermined configurations where $(\dim(W),\dim(Y))$ is $(4,4)$, $(16,16)$ or $(16,4)$,
with equally spaced measurement sites. In the underdetermined configurations $(4,16)$, the state and parameter error 
stagnates, while the residual error decreases to zero, which reflects the fact
that there are several $(y,u)\in Y\times w+W^\perp$ satisfying $\cR(y,u)=0$,
making the fundamental barrier $\delta_0$ strictly positive.
\item
The state error, residual and parameter error do not decrease to zero in the overdetermined configuration $(\dim(W),\dim(Y))=(8,4)$ where the measurement
sites are concentrated on the upper-half of the domain. This case is interesting
since, while we may expect that there is a unique pair $(y,u(y))\in Y\times w+W^\perp$
reaching the global minimal value $\cR(y,u)=0$, the algorithm seems to get trapped in local minima.
\end{enumerate}

\begin{figure}[H]
\centering   
  \includegraphics[width=\linewidth]{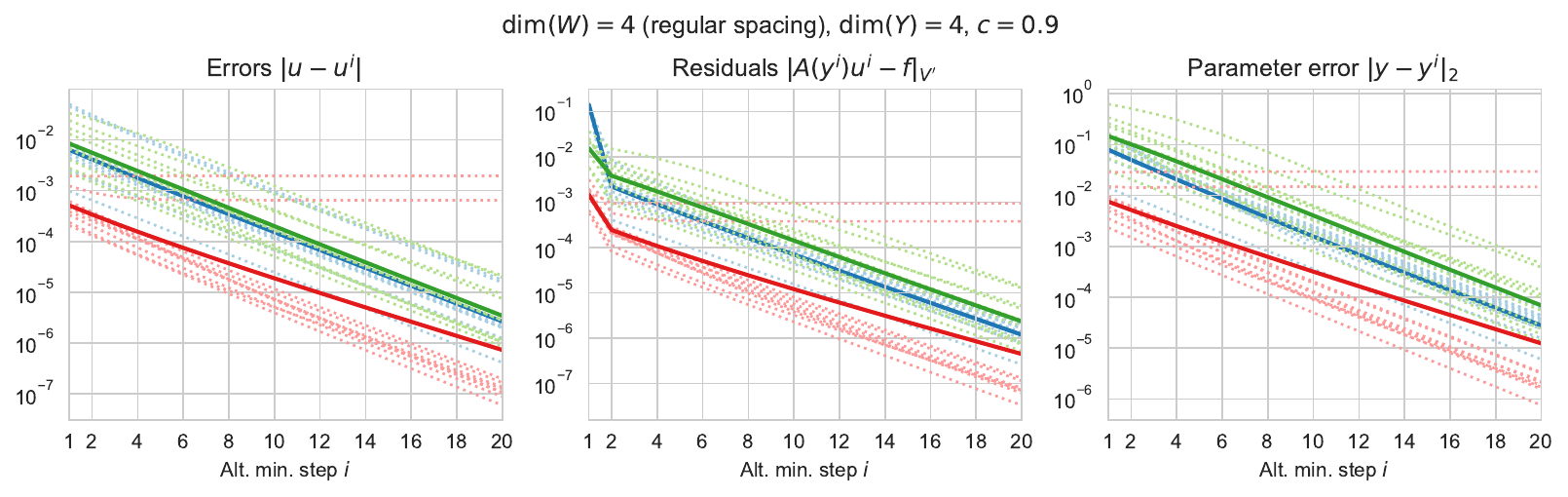}
  \includegraphics[width=\linewidth]{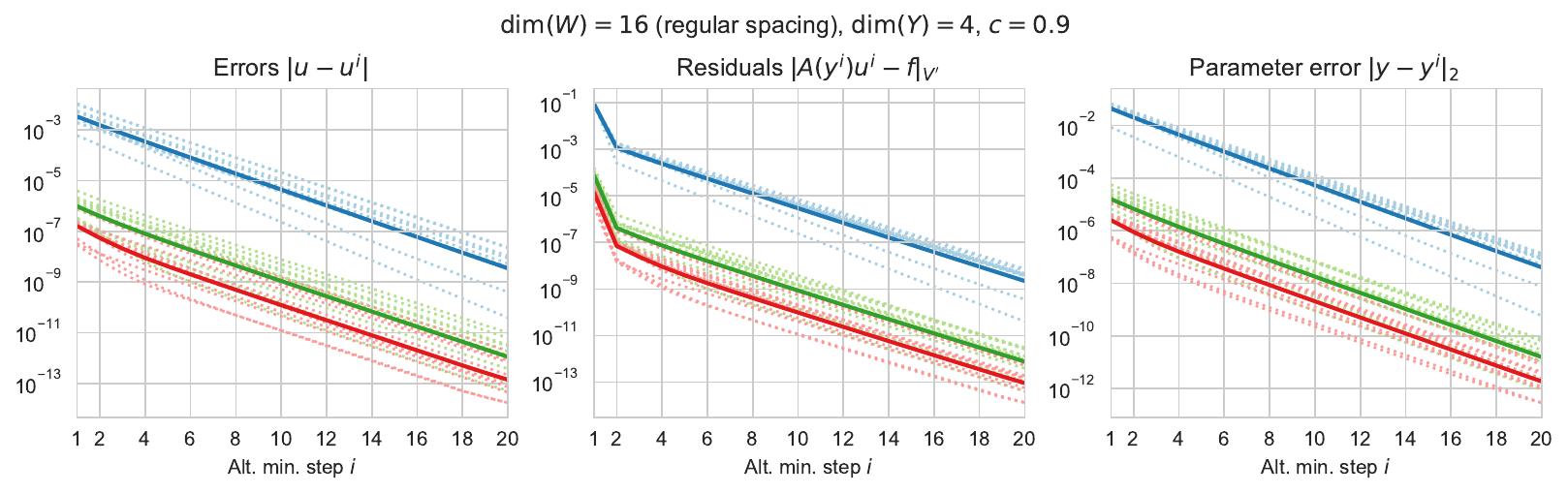}
  \includegraphics[width=\linewidth]{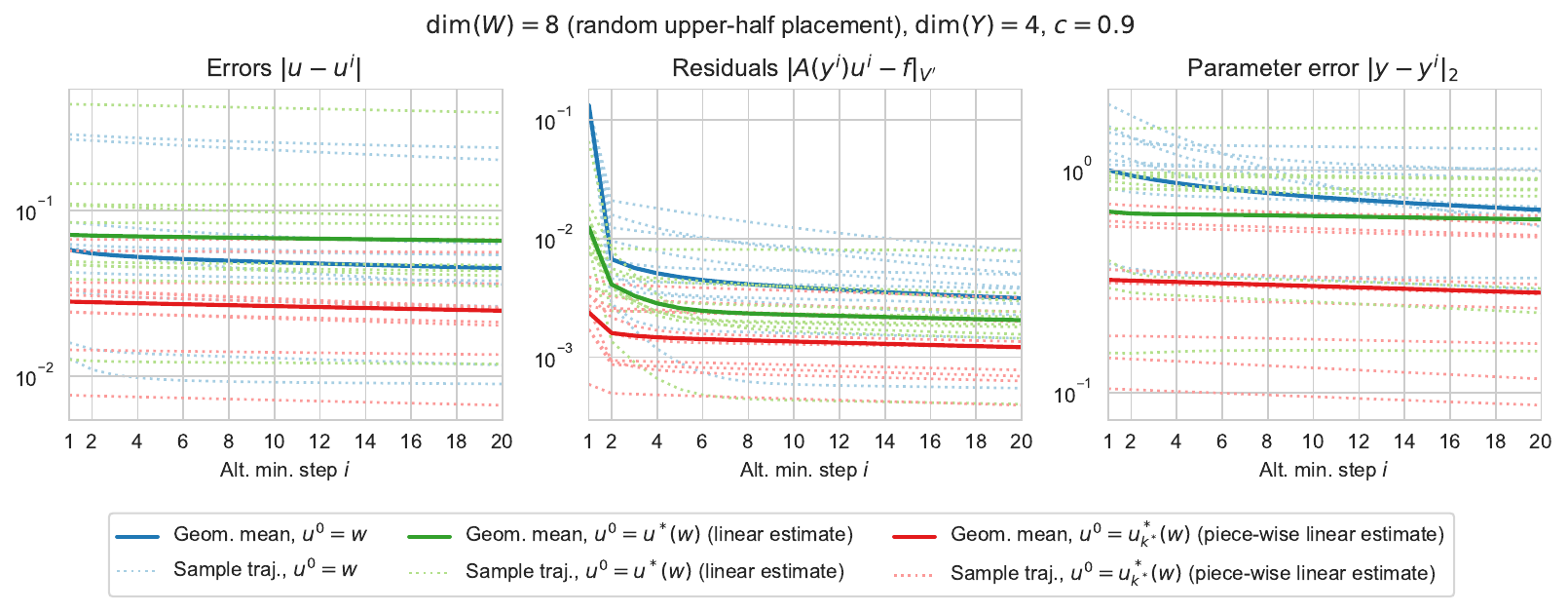}
\end{figure}

\begin{figure}
\centering   
  \includegraphics[width=\linewidth]{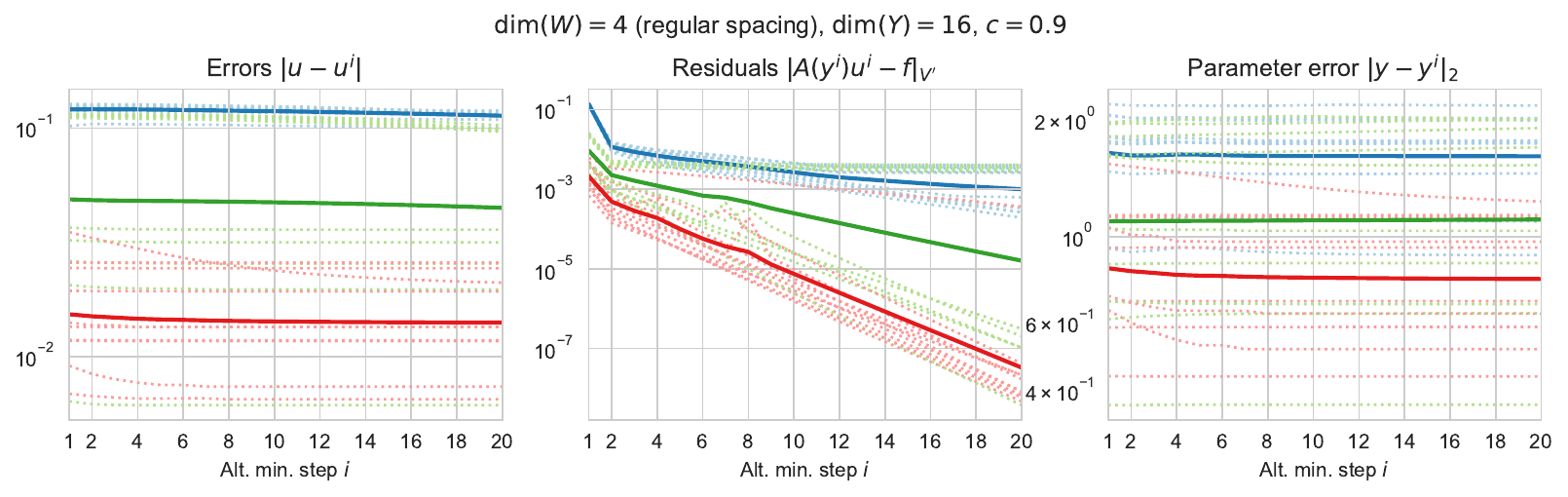}
  \includegraphics[width=\linewidth]{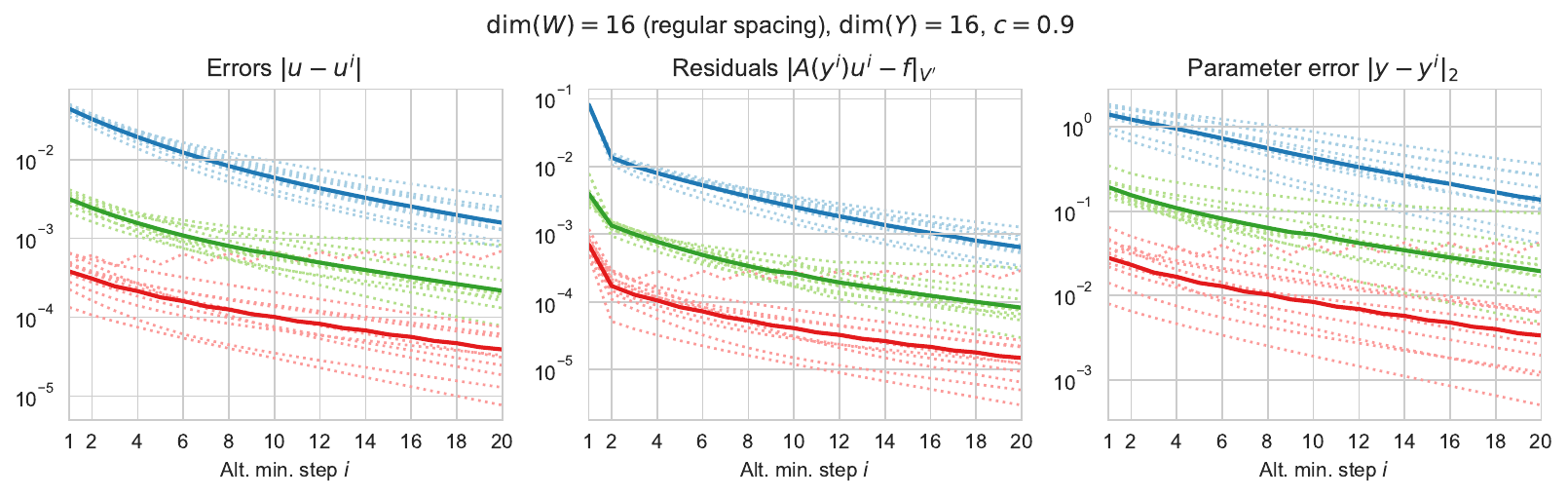}
  \includegraphics[width=\linewidth]{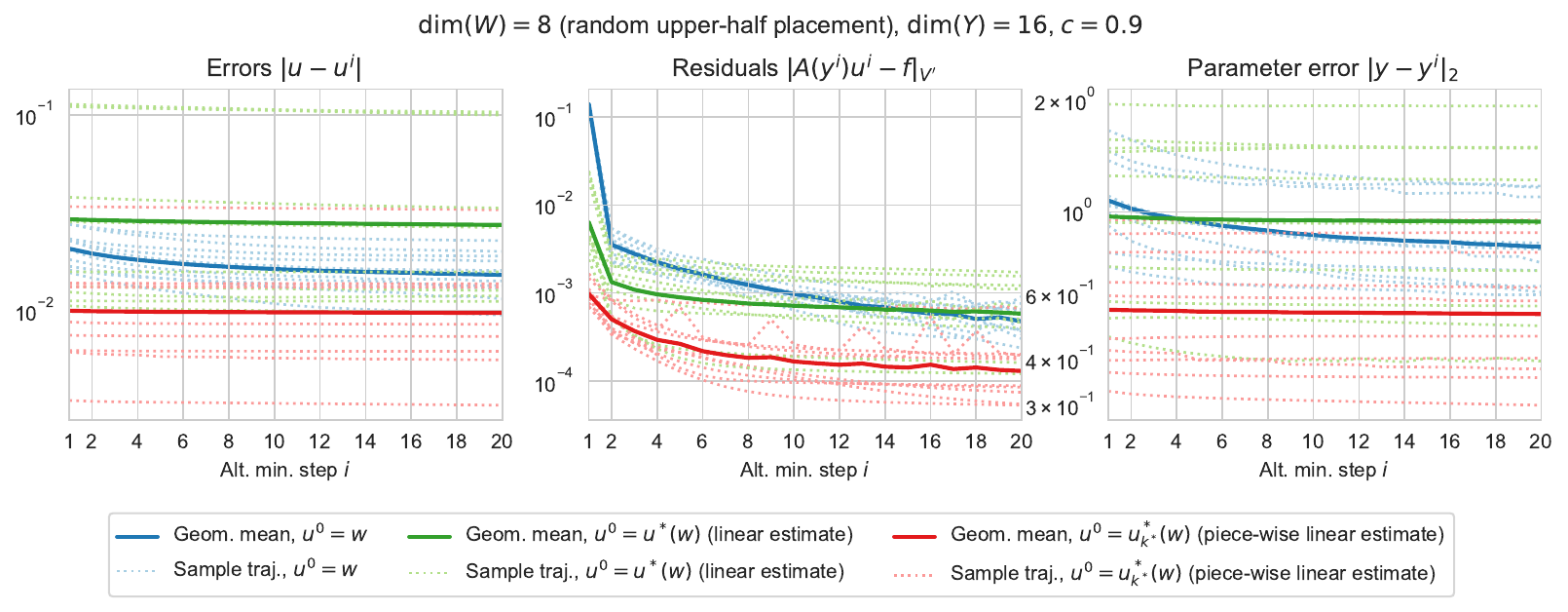}
\end{figure}

\end{document}